\documentclass[10pt]{amsart}

\usepackage{amsaddr}
\usepackage[dvips]{graphicx}
\usepackage{epsfig}
\usepackage{amsmath,amssymb,amscd,amsthm}

\usepackage{graphics}
\usepackage{picins}
\usepackage{colortbl}

\begin{document}

\setcounter{secnumdepth}{3}
\setcounter{tocdepth}{2}
\newtheorem{theorem}{Theorem}[section]
\newtheorem{definition}[theorem]{Definition}
\newtheorem{lemma}[theorem]{Lemma}
\newtheorem{proposition}[theorem]{Proposition}
\newtheorem{corollary}[theorem]{Corollary}
\newtheorem{gstatement}{Problem}
\newtheorem{question}{Question}
\newcommand{\mf}{\mathfrak}
\newcommand{\mb}{\mathbb}
\newcommand{\ol}{\overline}
\newcommand{\la}{\langle}
\newcommand{\ra}{\rangle}

\newtheorem{Alphatheorem}{Theorem}
\renewcommand{\theAlphatheorem}{\Alph{Alphatheorem}} 

\newtheorem{Alphatheoremprime}{Theorem}
\renewcommand{\theAlphatheoremprime}{\theAlphatheorem\textquotesingle}

\newcommand{\EM}{\ensuremath}
\newcommand{\norm}[1]{\EM{\left\| #1 \right\|}}

\newcommand{\modul}[1]{\left| #1\right|}

\title[From spaces of polygons to spaces of polyhedra]{From spaces of polygons to spaces of polyhedra following Bavard, Ghys and Thurston}
\author{Fran\c{c}ois Fillastre}
\address{
University of Cergy-Pontoise, UMR~CNRS~8088, departement of mathematics,
\\ F-95000 Cergy-Pontoise \\
 FRANCE\\ 
{\tt francois.fillastre@u-cergy.fr} \\
{\tt http://www.u-cergy.fr/ffillast/}
}
\date{\today}
\thanks{The author was partially supported by Schweizerischer Nationalfonds 200020-113199/1.} 

\maketitle

\begin{abstract}
 After  work of W.~P.~Thurston, C.~Bavard and \'E.~Ghys constructed particular hyperbolic polyhedra from spaces of deformations of Euclidean polygons. We present this construction as a straightforward consequence of the theory of mixed-volumes.

The gluing of these polyhedra can be isometrically embedded into complex hyperbolic cone-manifolds constructed by Thurston from spaces of deformations of Euclidean polyhedra. It is then possible to deduce the metric structure of the spaces of polygons  embedded in  complex hyperbolic orbifolds  discovered by P.~Deligne and G.~D.~Mostow. 

\end{abstract}

In \cite{T} W.P.~Thurston described a natural complex hyperbolic structure on the space of  convex polytopes in  Euclidean 3-space with fixed cone-angles. Applying this construction to polygons, C.~Bavard and \'E.~Ghys pointed out in \cite{BG}  that spaces of convex Euclidean polygons with fixed angles are isometric to  particular hyperbolic polyhedra, called (truncated) orthoschemes. In Section~\ref{sec : 1} we get the Bavard--Ghys results by using the theory of mixed-area (mixed-volume for polygons). Along the way we obtain Proposition~\ref{prop: ortho=pol} which is new.
The use of the Alexandrov--Fenchel Theorem might appear artificial at this point (see the discussion after Theorem~\ref{thm: alex fenchel}), but mixed-area theory sheds light on  the relations between convex polygons and hyperbolic orthoschemes \emph{via} Napier cycles, see Subsection~\ref{subsec: napier}. Moreover, it is very natural as mixed--area is the polar form of the quadratic form studied in \cite{T,BG}. Above all,  it  indicates a way to generalize the Bavard--Ghys construction from spaces of polygons to spaces of polytopes of any dimension $d$. In the case  $d=3$, the construction is related to Thurston's, but is different. Further explanations will be given in a forthcoming paper \cite{AF}. Section~\ref{sec : 1} ends with a discussion of hyperbolic orthoschemes which are of Coxeter type, as it appears that the list given by Im Hof in \cite{ImHof1} is incomplete.

In Section~\ref{sec : 2} we glue some of these hyperbolic orthoschemes to get hyperbolic cone-manifolds. This can be seen as  hyperbolization of the space of configurations of weighted points on the circle. This has been done several times, especially in lower dimensions, but it seems that  the link with orthoschemes was never clearly established. Proposition~\ref{lem: calcul angle} is new as stated for all dimensions.

Section~\ref{sec: 3} describes a local parametrization of spaces of polyhedra, equivalent to that in \cite{T}, which points up a ``complex mixed-area''. We outline next the remainder of the construction of \cite{T} which allows one to recover in a simple way complex hyperbolic orbifolds  listed by Mostow (our interest will be for a sublist first established by Deligne and Mostow \cite{DM}).

 Finally in Section~\ref{sec: 4} we check that the spaces of polygons embed (isometrically)    into the spaces of polyhedra locally as real forms. This is a non-surprising and certainly well-known fact (see for example \cite{KM1}), which is  checked here with the parametrizations  we defined.  
We then easily derive  Theorem~\ref{thm: table}, which gives the metric structure of those sets of polygons  in the Deligne--Mostow orbifolds (are they manifolds, orbifolds, or just cone-manifolds?). 

\textbf{Acknowledgments.} 
 I was introduced to the subject by Jean--Marc Schlenker. The link with the mixed-volume theory comes from discussions with Ivan Izmestiev after he presented the content of \cite{BobenkoIzmestiev}. The existence of the Coxeter polyhedron represented in Figure~\ref{fig: tmarkinpol} was communicated to me by Anna Felikson and Pavel Tumarkin. I had  fruitful  discussions with them, as well as with Christophe Bavard and Ruth Kellerhals. 

A part of this work
was completed during my visits to  the research group ``Polyhedral
surfaces'' at TU Berlin, which I want to thank for its
hospitality.

The author wishes to thank the anonymous referees and Hans Rugh for their  comments.

\section{Spaces of polygons and hyperbolic orthoschemes}\label{sec : 1}

\subsection{Basics about hyperbolic polyhedra, Napier cycles}\label{subsec: napier}

 The \emph{signature} $(N,Z,P)$ of a symmetric bilinear form (or of a Hermitian form) is the triple 
constituted of its  $N$ negative eigenvalues, $Z$ zero eigenvalues and $P$ positive eigenvalues (with multiplicity).
We denote by $\mathbb{R}^{n,1}$ the Minkowski space of dimension $(n+1)$, that is $\mathbb{R}^{n+1}$ endowed with the bilinear form of signature $(1,0,n)$
$$\langle x,y \rangle_1 = -x_0y_0+x_1y_1+\ldots +x_ny_n.  $$
A vector $x$ of Minkowski space is said to be \emph{positive} if $\langle x,x\rangle_1 >0$ and \emph{non-positive} otherwise. The hyperbolic space of dimension $n$ is the following submanifold of $\mathbb{R}^{n,1}$ together with the induced metric: 
$$\mathbb{H}^n:=\{x\in \mathbb{R}^{n,1} \vert \langle x,x\rangle_1=-1, x_0>0 \} . $$
A convex polyhedron $P$ of $\mathbb{H}^n$ is the non-empty intersection of $\mathbb{H}^n$ with a convex polyhedral cone of  $\mathbb{R}^{n,1}$ with vertex at the origin. If two facets (i.e. codimension 1 faces)  of $P$  intersect in $\mathbb{H}^n$, their outward  normals in  $\mathbb{R}^{n,1}$ span a Riemannian plane and the angle between these two vectors is the exterior dihedral angle between the facets. The interior dihedral angle between the facets is $\pi$ minus the exterior dihedral angle. In this paper, the \emph{dihedral angle} is the interior dihedral angle.

Let us consider the central projection in $\mathbb{R}^{n,1}$ onto the hyperplane $\{x_0=1\}$. The image of the hyperbolic space under this projection is the interior of the unit ball of $\mathbb{R}^{n}$. It is endowed with the metric for which the projection is an isometry. In this model, known as the Klein projective model of the hyperbolic space, geodesics are straight lines. The unit sphere in this model is the boundary at infinity of the hyperbolic space. In this paper we call a \textit{convex generalized polyhedron} of the hyperbolic space (the hyperbolic part of) a convex polytope of $\mathbb{R}^{n}$ such that all its edges meet the interior of the unit ball (a \textit{polytope} is a compact polyhedron). A vertex lying outside the interior of the ball is called \textit{hyperideal}. It is \textit{ideal} if it lies on the unit sphere and \textit{strictly hyperideal} otherwise. A vertex in the interior of the unit ball is called \textit{finite}. A strictly hyperideal vertex corresponds to a positive vector of $\mathbb{R}^{n,1}$. The polyhedron is \emph{truncated} if we cut it along the hyperplanes orthogonal to its strictly hyperideal vertices. We get a new hyperbolic polyhedron with  new facets, one for each strictly hyperideal vertex $v$. Such a new facet has the property of being orthogonal to all the facets which had $v$ as a vertex. A hyperbolic convex generalized polyhedron with only finite vertices is compact, and it is of finite volume if it has only finite and ideal vertices.

A particularly important class of hyperbolic polyhedra (compact or of finite volume) is that of  \emph{Coxeter polyhedra}, whose dihedral angles are integer submultiples of $\pi$. This implies in particular that the polyhedron is simple (this means that $n$ facets  meet at each finite vertex). For more details about Coxeter polyhedra we refer to \cite{vinberg,vinberglivre}. Coxeter polyhedra are represented  by  \emph{Coxeter diagrams}. Each facet is represented by a node. If two facets intersect orthogonally the nodes are not joined. If the two facets intersect at an angle $\pi/k$, $k>2$, the nodes are joined by a line with a $k$ above it. If the facets intersect at infinity we put a $\infty$ on the line, and if the facets do not intersect the nodes are joined by a  dashed line.

Let us consider a set of  vectors  $e_k\in\mathbb{R}^{n,1}$, $k\in\mathbb{Z}$ modulo $n+3$, such that
\begin{itemize}
 \item $\langle e_{k-1},e_k \rangle_1 <0$ for all $k$;
\item $\langle e_{k},e_j \rangle_1 =0$ for  $2\leq \vert j-k\vert\leq n+1$.
\end{itemize}
Two sets of vectors $\{e_k\}$ and $\{f_k\}$ as above are considered to be  equivalent if, 
for each $k$, $e_k=\lambda_k f_k$ with $\lambda_k$ a positive scalar. The equivalence class is a \emph{Napier cycle}.

The definition above comes from \cite{ImHof1}. As noted in  remarks on pp. 526 and 531 of this reference, this definition is suggested by the ``Napier pole sequences'' introduced in \cite{DB}. A Napier pole sequence is a sequence of unit vectors in Euclidean space satisfying the same equations as above (with the usual scalar product instead of $\langle .,. \rangle_1$). It is shown in \cite[Lemma 5.2]{DB} that  the sequence is then periodic. In our definition  the periodicity is assumed but it is not hard to see that it  is implied by the other assumptions, following the lines of the proofs of \cite[Lemma 5.2]{DB} and  \cite[Proposition 1.2]{ImHof1}.

Along the way we get that a Napier cycle always contains  $n+1$ consecutive  positive vectors, which generate the whole  cycle, and   there are three types of Napier cycles (see \cite{ImHof1}):
\begin{itemize}
 \item \emph{type 1}: two adjacent vectors are non-positive; 
 \item \emph{type 2}: one vector is non-positive;
 \item \emph{type 3}: all vectors  are positive.
\end{itemize}

There are  corresponding polyhedra, bounded by 
 the hyperplanes orthogonal to the positive vectors of the Napier cycles. If a Napier cycle has non-positive vectors, then they correspond to vertices of the polyhedron, see \cite[2]{ImHof1}. 

\begin{itemize}
 \item The (ordered) set of  outward normals of an \emph{ordinary orthoscheme}  
generates a  Napier cycle of type 1. (If it is Coxeter,) Its Coxeter diagram is a linear chain with $n+1$ nodes.   
 \item The set of  outward normals of a \emph{simply-truncated orthoscheme} 
generates a  Napier cycle of type 2. Its Coxeter diagram is a linear chain with $n+2$ nodes.   
 \item The  set of  outward normals of a  \emph{doubly-truncated orthoscheme} is a Napier cycle of type 3. Its Coxeter diagram is a cycle with $n+3$ nodes.   
\end{itemize}
By abuse of language we will call  a polyhedron of one of the three types above an \emph{orthoscheme}.  
Usually, the  word orthoscheme designates  what we called  ordinary orthoscheme. See \cite{ImHof1} for more details about the terminology.

\textbf{Example.} In $\mathbb{H}^2$, an ordinary orthoscheme is a right triangle, a simply-truncated orthoscheme is a quadrilateral with three right angles  and  a doubly-truncated orthoscheme is a right-angled pentagon, see the figures in \cite{ImHof1}.

\textbf{Remark.} In \cite{DB}, to a Napier pole sequence there is an associated  ``Napier configuration'', which is a sequence of spherical orthoschemes. Geometrically this means that some vertices of a given orthoscheme will be considered as outward normals of another orthoscheme (we refer to \cite{DB} for more details). Analogous considerations in our case of Napier cycles would oblige us to consider larger class of polyhedra than hyperbolic ones. 
The terminology  ``Napier pole sequence'' comes from the fact that in the sphere of dimension 2,  relations  in a Napier configuration are Napier's rules (Napier is sometimes written Neper), see \cite[5]{DB}.

\subsection{Euclidean polygons and Napier cycles}

\begin{figure}[ht] \begin{center}
\input{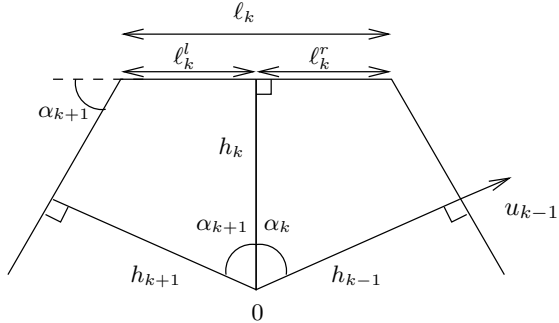}%
\caption{Notations for a convex polygon. \label{fig: def pol}}
\end{center}
\end{figure}

Let $P$ be a convex polygon of $\mathbb{R}^2$ with $n+3$ vertices such that the origin is contained in its interior.
We call $F_k$ the edges of $P$, labeled in  cyclic order, $\alpha_k$ is the exterior angle between $F_{k-1}$ and $F_k$, and $h_k$ is the distance of $F_k$ to the origin. The angles $\alpha_k$ satisfy
\begin{enumerate}
 \item[(A)] $0<\alpha_k<\pi$, $\displaystyle\sum_{k=1}^{n+3} \alpha_k=2\pi$.
\end{enumerate} 
We call $\ell_k$  the length of $F_k$, and we have (see Figure~\ref{fig: def pol})
\begin{equation}\label{eq:long}
\ell_k=\ell_k^r+\ell_k^l=\frac{h_{k-1}-h_k\cos(\alpha_{k})}{\sin(\alpha_{k})}+\frac{h_{k+1}-h_k\cos(\alpha_{k+1})}{\sin(\alpha_{k+1})}.
\end{equation}

We identify the set of \emph{heights} $h_1,\ldots,h_{n+3}$ (also called \emph{support numbers}) with $\mathbb{R}^{n+3}$ such that the set of outward unit normals $u_1,\ldots,u_{n+3}$ of the convex polygon $P$ corresponds to the canonical basis. For a vector $P\in\mathbb{R}^{n+3}$,   $h_k(P)$ is the 
$k$th coefficient of $P$ for the basis  $u_1,\ldots,u_{n+3}$. In particular  $h_k(u_i)=\delta^k_i$. We define $\ell_k(P)$ as the right-hand side of (\ref{eq:long}), after replacing the entries $h_i$ by $h_i(P)$. 

\begin{figure}[ht] \begin{center}
\input{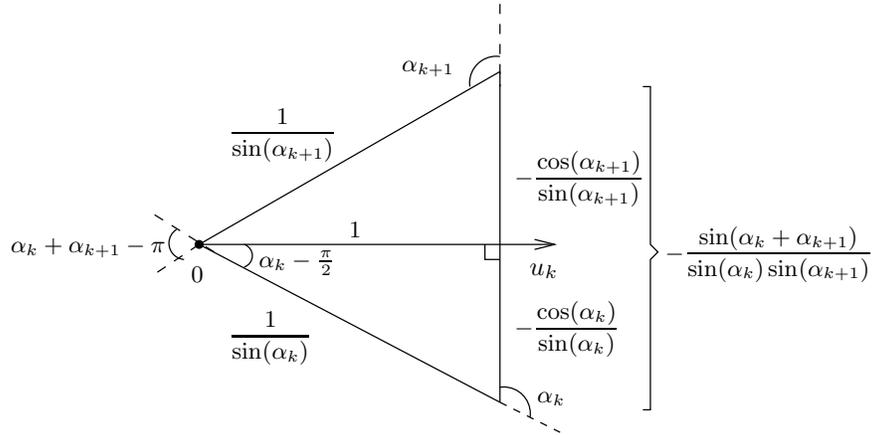}%
\caption{The geometric meaning of (\ref{eq:squarrednorm}): it is (minus) the signed area of the triangle described by $u_k$. \label{fig: triangle}}
\end{center}
\end{figure}

\textbf{A trivial example.} 
An element of $\mathbb{R}^{n+3}$ describes a polygon (not necessarily convex)
with $n+3$ edges (maybe of length $0$), which is such that the $k$th edge has outward normal $u_k$ and is on a line at distance $h_k$ from the origin. 
Let us consider the element  $u_k$ of $\mathbb{R}^{n+3}$. There is one edge on a line $l$ with normal $u_k$ and at distance $1$ from the origin. The edge with normal $u_{k-1}$ (resp. $u_{k+1}$) is at distance $0$ from the origin, hence it is on the unique line from the origin  making an angle $\alpha_k$ (resp. $\alpha_{k+1}$) with $l$. The other edges are reduced to the origin because they link the origin to itself. So  $u_k$ describes a triangle, see Figure~\ref{fig: triangle}.

We define the following bilinear form on $\mathbb{R}^{n+3}$:
$$m(P,Q):=-\frac{1}{2}\sum_{k=1}^{n+3}h_k(P)\ell_k(Q).$$

If $P$ and $Q$ are two convex polygons with outward unit normals $u_1,\ldots,u_{n+3}$, $m(P,Q)$ is known as (minus) the mixed-area of $P$ and $Q$, and $m(P,P)$ is minus the area of $P$, which is, up to the sign, the quadratic form used in \cite{BG} (in the present paper the minus sign serves only to get a more usual signature below). 
We get immediately that
\begin{equation}\label{eq:squarrednorm}
m(u_k,u_k)=-\frac{1}{2} \ell_k(u_k)=\frac{1}{2}\frac{\sin(\alpha_k+\alpha_{k+1})}{\sin(\alpha_k)\sin(\alpha_{k+1})}.\end{equation}
This formula has a geometric meaning: it is (minus) the signed area of the triangle described by $u_k$  (see Figure~\ref{fig: triangle} and the example above). 
 
It is also straightforward that 
\begin{equation}\label{eq: mui}
m(u_k,u_j)=\left\{
    \begin{array}{ccc}
        0 & \mbox{if } & 2\leq \vert j-k\vert\leq n+1 \\
      \displaystyle{  -\frac{1}{2}\frac{1}{\sin(\alpha_k)}} & \mbox{if} & j=k-1 \\
	 \displaystyle{-\frac{1}{2}\frac{1}{\sin(\alpha_{k+1})}}& \mbox{if} & j=k+1
    \end{array}
\right..
\end{equation}

It follows from (\ref{eq: mui}) that $m$ is symmetric, that is:
$$m(P,Q)=-\frac{1}{2}\sum_{k=1}^{n+3}h_k(Q)\ell_k(P).$$
This also follows from general properties of the mixed-volume \cite{schneider,AlexCP}.

\begin{theorem}\label{thm: alex fenchel}
 The symmetric bilinear form $m$ has signature $(1,2,n)$.
\end{theorem} 

Theorem~\ref{thm: alex fenchel} is a straightforward adaptation of the analogous result proved  in \cite{T} dealing with convex polytopes of $\mathbb{R}^3$ (see Section~\ref{sec: 3}). A more embracing statement is obtained in \cite{BG} with the same method  (see the remark after the proof of Proposition~\ref{prop: ortho=pol}). Theorem~\ref{thm: alex fenchel}  is also proved with greater effort in \cite{kojimaal1} in the case where all the $\alpha_k$ are equal. This statement is generalized to some cases of ``convex generalized polygons'' in \cite[Lemma 3.15]{BobenkoIzmestiev} (here ``generalized'' has another meaning than ours).

Theorem~\ref{thm: alex fenchel} is also a particular case of 
classical results about mixed-volumes, even if it is far from the simplest way of proving it. It appears in  Alexandrov's proof of the so-called  Alexandrov--Fenchel Theorem (or Alexandrov--Fenchel Inequality) for convex polytopes 
of $\mathbb{R}^d$ (here  $d=2$)  \cite{AlexAF,schneider,Alexcoll1}. Actually  Theorem~\ref{thm: alex fenchel} can be derived from 
 the Minkowski Inequality for convex polygons \cite[Note 1 p. 321]{schneider}, \cite{KlainMinkowski}: if $P$ and $Q$ are convex polygons then 
\begin{equation}\label{eq : mink}
 m(P,Q)^2\geq m(P,P)m(Q,Q)
\end{equation}
and equality occurs if and only if $P$ and $Q$ are homothetic. The way to go from Minkowski Inequality to Theorem~\ref{thm: alex fenchel} is a part of Alexandrov's proof of the Alexandrov--Fenchel Theorem. This is also done in a wider context in \cite[Appendix A]{Izmestiev}. Note that the Minkowski Inequality (\ref{eq : mink}) can be thought of as  the reversed Cauchy--Schwarz inequality in the case where the time-like vectors correspond to  convex polygons. One interest in considering  Bavard--Ghys construction as  a consequence of the Alexandrov--Fenchel Theorem in the particular case when $d=2$ rather than as a consequence of Thurston construction  is that 
 the former one generalizes immediately to any dimension $d$. This generalization will be the subject of \cite{AF}. 

The height of the sum of two polygons is the sum of the heights of the polygons, see e.g. \cite[Chapter IV]{Alexcoll1}. It follows that a translation of a polygon $P$ is the same as adding to the heights of $P$ the heights of a point. Hence the kernel of $m$ consists of heights spanning a point, since the area is invariant under translations (in the general case this is a step in the proof of the Alexandrov--Fenchel Theorem, see \cite[Lemma III p.71 ]{Alexcoll1}, \cite[Proposition 3 p.329]{schneider}). 

\begin{figure}[ht] \begin{center}
\input{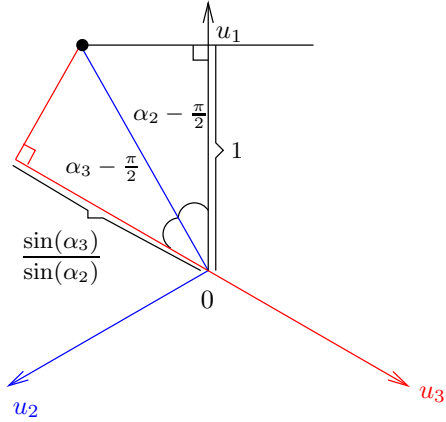}%
\caption{ The polygon with heights $1$, $0$ and  $-\frac{\sin(\alpha_3)}{\sin(\alpha_2)}$ is   a point. \label{fig: pointtriangle}}
\end{center}
\end{figure}

\textbf{A trivial example.} Let us consider $\mathbb{R}^3$, with basis $\{u_1,u_2,u_3\}$. From (\ref{eq:squarrednorm}), (\ref{eq: mui}) and (A) we get easily that the matrix of $m$ is
$$-\frac{1}{2}  \begin{pmatrix}
\frac{\sin(\alpha_3)}{\sin(\alpha_1)\sin(\alpha_2)} & \frac{1}{\sin(\alpha_2)}& \frac{1}{\sin(\alpha_1)}\\
\frac{1}{\sin(\alpha_2)}&\frac{\sin(\alpha_1)}{\sin(\alpha_2)\sin(\alpha_3)} & \frac{1}{\sin(\alpha_3)} \\
\frac{1}{\sin(\alpha_1)} & \frac{1}{\sin(\alpha_3)} & \frac{\sin(\alpha_2)}{\sin(\alpha_3)\sin(\alpha_1)}
\end{pmatrix} $$
and  the vectors

$$v:=\begin{pmatrix} 1 \\ 0 \\ -\frac{\sin(\alpha_3)}{\sin(\alpha_2)} \end{pmatrix}, \begin{pmatrix} 0 \\ 1 \\ -\frac{\sin(\alpha_1)}{\sin(\alpha_2)} \end{pmatrix}$$
span the kernel of $m$.
In Figure~\ref{fig: pointtriangle} we check that the polygon corresponding to the vector $v$ is a single point.
Another eigenvector of $m$ is
$$\begin{pmatrix} 1 \\ \frac{\sin(\alpha_1)}{\sin(\alpha_3)} \\ \frac{\sin(\alpha_2)}{\sin(\alpha_3)} \end{pmatrix}$$
 with negative eigenvalue
$$ -\frac{1}{2}\frac{\sin(\alpha_1)^2+\sin(\alpha_2)^2+\sin(\alpha_3)^2}{\sin(\alpha_1)\sin(\alpha_2)\sin(\alpha_3)},$$
so the signature of $m$ is $(1,2,0)$.

By Theorem~\ref{thm: alex fenchel}  the quotient of $\mathbb{R}^{n+3}$ by the kernel of $m$ is isometric to the Minkowski space $\mathbb{R}^{n,1}$. We denote by $\Pi$ the quotient map.
By definition it follows from (\ref{eq: mui}) and Theorem~\ref{thm: alex fenchel} that:
\begin{corollary}\label{prop: pol=napier}
 The set of vectors $(\Pi(u_1),\ldots,\Pi(u_{n+3}))$ is a Napier cycle.
\end{corollary}

Up to global isometries, this Napier cycle depends only on the angles $\alpha_k$ between the $u_k$'s. 
Let   $(\alpha_1,\ldots,\alpha_{n+3})$ be an ordered list of real numbers, up to cyclic permutations,  satisfying (A). 
We denote by $\mathcal{H}(\alpha_1,\ldots,\alpha_{n+3})$ the hyperbolic orthoscheme corresponding to  the Napier cycle given by the corollary above. We check now  that, as announced in the title, $\mathcal{H}(\alpha_1,\ldots,\alpha_{n+3})$ is a space of convex polygons.

\begin{lemma}\label{lem: ortho=pol}
 The orthoscheme $\mathcal{H}(\alpha_1,\ldots,\alpha_{n+3})$ is in bijection with the set of convex polygons with the $k$th exterior angle equal to  $\alpha_k$, up to direct isometries and homotheties.
\end{lemma}

 \begin{proof}
The interior of the hyperbolic polyhedron $\mathcal{H}(\alpha_1,\ldots,\alpha_{n+3})$   is the set of vectors $P$ satisfying $m(P,u_k)<0$ for all positive $u_k$, i.e. the set of $P$ such that $\ell_k(P)>0$ (it is not hard to see that this is true even if one or two $u_k$ are non-positive). So $P$ belongs to  the set of convex polygons with the $u_k$ as outward unit normals, which is  the set of convex polygons with $k$th exterior angle equal to  $\alpha_k$, up to rotations as the $u_k$ are arbitrarily placed in the plane. We saw that to quotient by the kernel of  $m$ is equivalent to consider the polygons up to translation. 
Finally  $\mathcal{H}(\alpha_1,\ldots,\alpha_{n+3})$ is a subset of the hyperbolic space, hence it contains only polygons of unit area, that is the same as polygons up to homotheties.
\end{proof}
From now we identify the space of polygons described in the lemma with  $\mathcal{H}(\alpha_1,\ldots,\alpha_{n+3})$.

\begin{figure}[ht] \begin{center}
\input{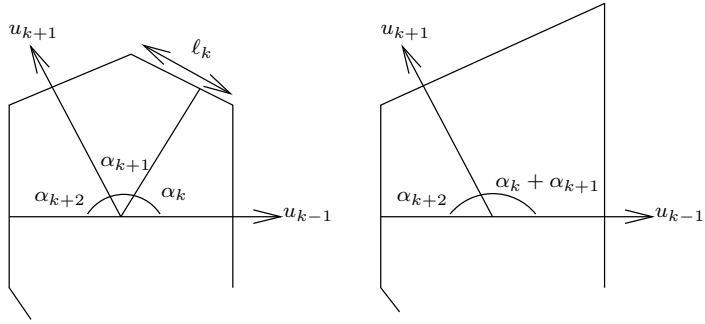}%
\caption{ On the left $u_{k+1}$ is a space-like vector and it becomes a light-like vector on the right because $\alpha_k+\alpha_{k+1}+\alpha_{k+2}=\pi$.
\label{fig: lightlike}}
\end{center}
\end{figure}

\begin{proposition}\label{prop: ortho}
Suppose $n\geq 2$. The hyperbolic orthoscheme $\mathcal{H}(\alpha_1,\ldots,\alpha_{n+3})$ has the following properties.
\begin{itemize}
\item[i)] Its type is equal to $3$ minus the number of $k$ such that $\alpha_k+\alpha_{k+1}\geq\pi$.
  \item[ii)] It has non-obtuse dihedral angles, and if  $\Pi(u_{k-1})$ and $\Pi(u_{k})$ are space-like, the dihedral angle $\Theta$  between the corresponding facets is acute and satisfies
\begin{equation}\label{eq: cos angle} \cos^2(\Theta)=\frac{\sin(\alpha_{k-1})\sin(\alpha_{k+1})}{\sin(\alpha_{k-1}+\alpha_{k})\sin(\alpha_{k}+\alpha_{k+1})}.\end{equation}
\item[iii)] 
It is of finite volume. Moreover  it is compact if and only if there is no couple  $k,k'$ for which $\alpha_k+\cdots+\alpha_{k'}=\pi$. \end{itemize}
\end{proposition}
The last condition about compactness can be rephrased by saying that the polygons of $\mathcal{H}(\alpha_1,\ldots,\alpha_{n+3})$ have no parallel edges.
The cases $n=0$ and $n=1$ are geometrically	 meaningless. Note that if $\alpha,\beta,\gamma,\delta$ are the angles of a parallelogram, then $\mathcal{H}(\alpha,\beta,\gamma,\delta)$ is the whole ``hyperbolic line'' $\mathbb{H}^1$.
\begin{proof}
First note that due to (\ref{eq:squarrednorm}) the character of the vector $\Pi(u_k)$ is easily characterized:
\begin{itemize}
 \item it is space-like if $\alpha_k+\alpha_{k+1}<\pi$,
\item it is light-like if $\alpha_k+\alpha_{k+1}=\pi$,
\item it is time-like if $\alpha_k+\alpha_{k+1}>\pi$.
\end{itemize}
and i) follows. The dihedral angles are either $\pi/2$ else  minus the cosine of the angle is given by 
\begin{equation}\label{eq: normalisation}
\frac{m(u_{k-1},u_{k})}{\sqrt{m(u_{k-1},u_{k-1})}\sqrt{m(u_{k},u_{k})}}=-\sqrt{\frac{\sin(\alpha_{k-1})\sin(\alpha_{k+1})}{\sin(\alpha_{k-1}+\alpha_{k})\sin(\alpha_{k}+\alpha_{k+1})}},\end{equation}
 which is a real negative number if $\alpha_{k-1}+\alpha_{k}<\pi$ and $\alpha_{k}+\alpha_{k+1}<\pi$, that proves ii). Actually   orthoschemes have always non-obtuse dihedral angles and finite volume  \cite[Proposition 2.1]{ImHof1}. This gives the first part of iii).
It remains to prove the assertion about the compactness. The polyhedron is not compact if and only if one of its faces contains an ideal point. It can be the case if a vector $u_k$ is light-like, that is  if  $\alpha_k+\alpha_{k+1}=\pi$.  Suppose now that $u_k$ is space-like. A facet  of the polyhedron is given by the polygons $P$ which satisfy $\ell_k(P)=0$. This facet is itself a polyhedron of one dimension lower. More precisely it is isometric to $\mathcal{H}':=\mathcal{H}(\alpha_1,\ldots,\alpha_{k-1},\alpha_{k}+\alpha_{k+1},\alpha_{k+2},\ldots,\alpha_{n+3})$. The hyperbolic polyhedron $\mathcal{H}'$ corresponds to a cone in an  ambient Minkowski space, in which (the vector corresponding to) $u_{k+1}$  is  light-like  if and only if $\alpha_k+\alpha_{k+1}+\alpha_{k+2}=\pi$. This is  easy to check from the definition of the bilinear form,  see also Figure~\ref{fig: lightlike}. Now if $u_{k+1}$ is space-like in $\mathcal{H}'$, we repeat the reasoning  for a facet of $\mathcal{H}'$, and so on. 
\end{proof}
We  easily check  below that Equation~(\ref{eq: cos angle}) is another form of the main formula  of \cite{BG} (written in the present paper as Equation~(\ref{eq: crossrat})). It also appears in the form of (\ref{eq: cos angle}) in \cite{kojimaal1,MorinHishi} for $n=2,3$.  We denote by $U_k$ the line spanned by the vector $u_k$ in $\mathbb{R}^2$, and we define the cross-ration as
$$[a,b,c,d]=\frac{d-a}{a-b}\frac{b-c}{c-d}. $$

\begin{corollary}
 Let $\Theta$ be as in \emph{ii)} of Proposition~\ref{prop: ortho}. Then
\begin{equation}\label{eq: crossrat}
\tan^2(\Theta)=-[U_{k-1},U_{k},U_{k+1},U_{k+2}].
\end{equation}
\end{corollary}
\begin{proof}
 This follows from (\ref{eq: cos angle}) and from the well-known fact that
\begin{equation}\label{eq: proj formula}
\frac{\sin(\alpha_{k-1})\sin(\alpha_{k+1})}{\sin(\alpha_{k-1}+\alpha_{k})\sin(\alpha_{k}+\alpha_{k+1})}=\frac{1}{1-[U_{k-1},U_{k},U_{k+1},U_{k+2}]}.
 \end{equation}
\end{proof}
This formula is stated in \cite{BG} for the directions of the edges of the polygon,  and not for the  directions of the normals of the polygon as above, but the cross-ratios are the same. 

We proved, following \cite{BG}, that from any convex polygon one can construct a hyperbolic orthoscheme. 
It is possible to prove the converse. This point is missing from \cite{BG}, but this reference contains the main point for the proof. Let us be more formal.

 Let $A(n)$ be the space of all the ordered lists $(\alpha_1,\ldots,\alpha_{n+3})$ up to the action of the dihedral group $\mathbb{D}_{n+3}$,
with $\alpha_k$  real numbers satisfying (A). Note that if $\sigma\in\mathbb{D}_{n+3}$  then  $\mathcal{H}(\alpha_1,\ldots,\alpha_{n+3})$ and $\mathcal{H}(\alpha_{\sigma(1)},\ldots,\alpha_{\sigma({n+3})})$ are isometric. We say that two elements $(\alpha_1,\ldots,\alpha_{n+3})$ and $(\alpha_1',\ldots,\alpha_{n+3}')$  of $A(n)$ are equivalent if, given  a set of planar vectors  $\mathcal{u}=\{u_1,\ldots,u_{n+3}\}$ with $\alpha_k$ the angle between $u_{k-1}$  and $u_k$ and  a set of planar vectors  $\mathcal{u}'=\{u_1',\ldots,u_{n+3}'\}$ with  $\alpha_k'$ the angle between $u_{k-1}'$  and $u_k'$, then  there exists a projective map $\varphi$ of $\mathbb{R}^2$ fixing the origin and sending $\mathcal{u}$ to $\mathcal{u}'$, i.e. $\varphi(u_k)=u'_k$ for all $k$. 
We denote by $\overline{A}(n)$ the quotient of $A(n)$ by this equivalence relation, and by $H(n)$ the space of orthoschemes of $\mathbb{H}^n$ up to global isometries.

\begin{proposition}\label{prop: ortho=pol}
There is a bijection between $\overline{A}(n)$ and $H(n)$.
\end{proposition}

\begin{proof}
Let $(\alpha_1,\ldots,\alpha_{n+3})\in A(n)$. We know that these numbers define a Napier cycle. The corresponding orthoscheme is defined by the hyperplanes orthogonal to the positive vectors of the Napier cycle. Hence it is defined by the Gram matrix of the normalized  positive vectors of the Napier cycle, whose coefficients are either $1$ or given by (\ref{eq: normalisation}), which is a projectively invariant formula (see (\ref{eq: proj formula})). Hence there is a well-defined map from 
$\overline{A}(n)$ to $H(n)$. 

Let $H\in H(n)$ and let  $(e_1,\ldots,e_{n+3})$ be the  Napier cycle generated by the outward unit normals of $H$. We know that there are at most two non-positive vectors. Moreover if there are two non-positive vectors, they are consecutive. In this case, up to change the labeling, we suppose that the non-positive vectors are $e_{n+2}$ and $e_{n+3}$. If there is only one non-positive vector, we suppose that it is $e_{n+3}$.   Suppose now that the two elements $(\alpha_1,\ldots,\alpha_{n+3})$ and $(\alpha_1',\ldots,\alpha_{n+3}')$ of $A(n)$  lead to $H$. Up to a projective transformation we can consider that $\alpha_1=\alpha_1'$ and $\alpha_2=\alpha_2'$. We denote by $s(.,.,.)$ the function defined by the right side of (\ref{eq: normalisation}). As $s(\alpha_1,\alpha_2,\alpha_3)$ and $s(\alpha_1',\alpha_2',\alpha_3')$ are both equal to $m(e_1,e_2)$ it follows easily that $\alpha_3=\alpha_3'$ (this is straightforward with the cross-ratio). 
Next, as $s(\alpha_2,\alpha_3,\alpha_4)$ and $s(\alpha_2',\alpha_3',\alpha_4')$ are both equal to $m(e_2,e_3)$ it follows  that $\alpha_4=\alpha_4'$, and so on until $s(\alpha_n,\alpha_{n+1},\alpha_{n+2})=m(e_{n},e_{n+1})=s(\alpha_n,\alpha_{n+1},\alpha_{n+2}')$, which gives $\alpha_{n+2}=\alpha_{n+2}'$. The last equality $\alpha_{n+3}=\alpha_{n+3}'$ follows as the sum of the angles is equal to $2\pi$. Hence the map from $\overline{A}(n)$ to $H(n)$ is injective.

For an arbitrary $H\in H(n)$, we can find  angles $(\alpha_1,\ldots,\alpha_{n+3})$ in a similar way: $\alpha_1,\alpha_2$ and $\alpha_3$ are chosen such that they are between $0$ and $\pi$ and such that $s(\alpha_1,\alpha_2,\alpha_3)=m(e_1,e_2)$ (it is clear that such a triple always exists). It is easy to see that the angle $\alpha_4$ (supposed to be between $0$ and $\pi$) is determined by  $s(\alpha_2,\alpha_3,\alpha_4)=m(e_2,e_3)$  and so on until $s(\alpha_n,\alpha_{n+1},\alpha_{n+2})=m(e_{n},e_{n+1})$ which determines $0<\alpha_{n+2}<\pi$.
Now we have to examine different cases, according to the character of $e_{n+2}$:
\begin{itemize}
 \item if $e_{n+2}$ is space-like, $0<\alpha_{n+3}<\pi$ is given as above;
\item if $e_{n+2}$ is light-like, we define $\alpha_{n+3}=\pi-\alpha_{n+2}$. It is between $0$ and $\pi$;
\item if $e_{n+2}$ is time-like, we again define $\alpha_{n+3}$ with the help of  (\ref{eq: normalisation}). The squared norm of $e_{n+1}$ is positive and the one of $e_{n+2}$ is negative. As they are elements of a Napier cycle, $m(e_{n+1},e_{n+2})$ is negative, hence the right-hand side of (\ref{eq: normalisation}) must belong to $i\mathbb{R}_-$.  As $\sin(\alpha_{n+1})$ and $\sin(\alpha_{n+1}+\alpha_{n+2})$ are positive and as $\sin(\alpha_{n+2}+\alpha_{n+3})$ is negative (because $e_{n+2}$ is time-like),  $\sin(\alpha_{n+3})$ must be positive, so that $0<\alpha_{n+3}<\pi$. 
\end{itemize}

Now we have $n+3$ angles between $0$ and $\pi$. For each such set of angles, (\ref{eq:squarrednorm})  and (\ref{eq: mui}) allow to define a bilinear form   (geometrically the signed area of the polygons constructed of the angles). In our case this form has signature $(1,2,n)$, because it is the Gram matrix of a Napier cycle. But it is known that such forms have this signature only if the $\alpha_i$ sum up to $2\pi$ \cite[Proposition p. 209]{BG}. Hence we constructed an element of $A(n)$.
\end{proof}

\textbf{Remark.} It could be interesting to know if some of the numerous results about hyperbolic orthoschemes can be translated in terms of Euclidean polygons (for such results, see \cite{DB} and references therein). Moreover  \cite{BG} contains also a computation of the signature of the (signed) area form for spaces of non-convex polygons. In particular one can construct Euclidean and spherical  polyhedra.

\subsection{Coxeter orthoschemes}

The aim of \cite{ImHof1} (previously announced in \cite{ImHof2}) is to find all the Coxeter orthoschemes. 
Some subfamilies were already known, especially in dimensions  $2$ and $3$ (see \cite{ImHof2,ImHof1} and \cite{vinberg2,vinberg} for more details). For these dimensions there exists infinite families of Coxeter orthoschemes. For dimension $\geq 4$, Im Hof found a list of $75$ Coxeter orthoschemes up to dimension $9$ (he proved that they can't exist for higher dimension). In \cite{BG} the existence of Coxeter orthoschemes is checked, by showing for each one a list of real numbers representing the slopes of the lines parallel to the edges of a suitable convex polygon.

\begin{figure}[ht] \begin{center}
\input{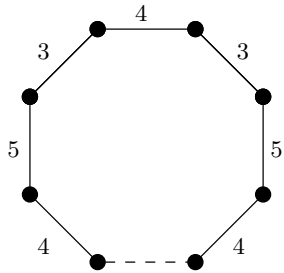}%
\caption{The Tumarkin polyhedron is a compact Coxeter orthoscheme of dimension $5$ and type $3$. \label{fig: tmarkinpol}}
\end{center}
\end{figure}

Figure~\ref{fig: tmarkinpol} represents another example of Coxeter orthoscheme coming from \cite{tumarkin}. We  name it \emph{Tumarkin polyhedron}. It is  not hard to find a convex polygon $P$ of $\mathbb{R}^2$ whose  space of angle-preserving deformations is isometric to the  Tumarkin polyhedron. The following list  gives the slopes of the lines containing the normals $u_k$ of $P$:
$$ \left(\sqrt{5},-2,-1,0,1,\infty,-3,\frac{\sqrt{5}-3}{2}\right).$$
This list confirms the existence of  Tumarkin polyhedron (dihedral angles can be easily computed with (\ref{eq: crossrat}). The fact that all the slopes are different indicates that the polyhedron is compact). This polyhedron doesn't appear in \cite{ImHof1} (nor in \cite{BG}). It seems that it has been forgotten
by Im Hof, and then it is natural to ask:
\begin{question}
 Is Im Hof's list together with the Tumarkin polyhedron complete?
\end{question}

\begin{figure}[ht] \begin{center}
\input{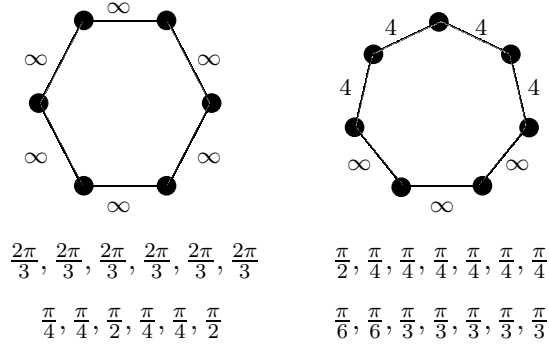}%
\caption{The Coxeter diagram on the left represents a Coxeter orthoscheme of dimension $3$ and type $3$ (the one shown in \cite[Figure 3]{T}). The  Coxeter diagram on the right represents a Coxeter orthoscheme of dimension  $4$ and type $3$. Below each diagram we give two different lists of rational angles $(\alpha_1,\ldots,\alpha_{n+3})$ such that  $\mathcal{H}(\alpha_1,\ldots,\alpha_{n+3})$ is isometric to the orthoscheme. \label{fig: 2ex cox}}
\end{center}
\end{figure}

We call a \emph{rational angle} an angle of the form $q\pi$, $q\in\mathbb{Q}$.
If $H$ is an orthoscheme, we can't hope that there exists a set of rational angles  $(\alpha_1,\ldots,\alpha_n)$ such that $\mathcal{H}(\alpha_1,\ldots,\alpha_{n+3})$ is isometric to $H$.
A natural question is if $H$ is a Coxeter orthoscheme, but I failed to find such set of angles for the Tumarkin polyhedron. Conversely two different sets of rational angles can lead to the same Coxeter orthoscheme,  examples are shown in Figure~\ref{fig: 2ex cox}.


\section{Spaces of polygons and hyperbolic cone-manifolds}\label{sec : 2}

\begin{figure}[ht] \begin{center}
\input{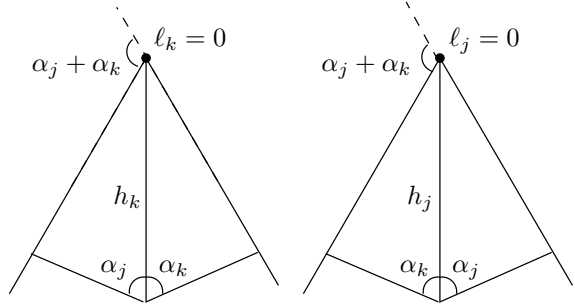}%
\caption{
Let $\alpha_k+\alpha_j<\pi$. A facet of $\mathcal{H}(\alpha_1,\ldots,\alpha_k,\alpha_j,\ldots,\alpha_n)$ is defined by $m(u_k,P)=0$, that is $\ell_k(P)=0$ and this space is isometric to 
$\mathcal{H}(\alpha_1,\ldots,\alpha_k+\alpha_j,\ldots,\alpha_n)$.
This last space is also isometric to the facet of 
 $\mathcal{H}(\alpha_1,\ldots,\alpha_j,\alpha_k,\ldots,\alpha_n)$ defined by $m(u_j,P)=0$.
  \label{fig: facet}}
\end{center}
\end{figure}

In the whole section we suppose that $n\geq2$.
Let $\alpha_1,\ldots,\alpha_{n+3}$ be $n+3$ real numbers  satisfying (A). 
We denote by $R(\alpha_1,\ldots,\alpha_{n+3})$  the set of all permutations of the set $\{\alpha_1,\ldots,\alpha_{n+3}\}$  up to the action of the dihedral group. 
Hence $R(\alpha_1,\ldots,\alpha_{n+3})$ has $(n+2)!/2$ elements. 
Each element  $\sigma\in R(\alpha_1,\ldots,\alpha_{n+3})$ will be identified with the orthoscheme $\mathcal{H}(\alpha_{\sigma(1)},\ldots,\alpha_{\sigma(n+3)})$.  There is a natural way to glue all these polyhedra to each other.
Let $\alpha_k,\alpha_j$ be such that $\alpha_k+\alpha_j<\pi$ (such pair always exists as the $\alpha_k$ satisfy (A)). It is easy to see that both orthoschemes $\mathcal{H}(\alpha_1,\ldots,\alpha_k,\alpha_j,\ldots,\alpha_n)$ and $\mathcal{H}(\alpha_1,\ldots,\alpha_j,\alpha_k,\ldots,\alpha_n)$  have a facet isometric to $\mathcal{H}(\alpha_1,\ldots,\alpha_k+\alpha_j,\ldots,\alpha_n)$, see Figure~\ref{fig: facet}. We glue them isometrically along this facet. We do so for all facets of all orthoschemes in  $R(\alpha_1,\ldots,\alpha_{n+3})$. At the end we obtain a space which is by construction a hyperbolic cone-manifold of dimension $n$. We still denote by $R(\alpha_1,\ldots,\alpha_{n+3})$  this space.

\begin{lemma}
 The cone-manifold  $R(\alpha_1,\ldots,\alpha_{n+3})$ is connected. It has finite volume and it is compact if and only if there is no couple $k,k'$ for which $\alpha_k+\cdots+\alpha_{k'}=\pi$.
\end{lemma}
Finiteness of the volume and description of the compactness are straightforward consequences of Proposition~\ref{prop: ortho}. In order to prove  connectedness
 we will prove the following lemma, which implies the one above.  We denote by $\tilde{R}(\alpha_1,\ldots,\alpha_{n+3})$ the double-covering of $R(\alpha_1,\ldots,\alpha_{n+3})$ obtained by distinguishing a list from the one obtained  by reversing the order.
\begin{lemma}\label{lem: connexite}
 The cone-manifold $\tilde{R}(\alpha_1,\ldots,\alpha_{n+3})$ is not connected if and only if  there exists $\alpha_i,\alpha_j,\alpha_k$ ($i\not= j,j\not=k,k\not= i$) such that
$$\alpha_i+\alpha_j\geq\pi, \alpha_j+\alpha_k\geq\pi,\alpha_k+\alpha_i\geq\pi. $$
 In this case $\tilde{R}(\alpha_1,\ldots,\alpha_{n+3})$ has two connected components, which are identified by reversing the order of the angles.
\end{lemma}
Note that this last case can happen only if $\tilde{R}(\alpha_1,\ldots,\alpha_{n+3})$ contains orthoschemes of type 1.
\begin{proof}
  To prove the lemma it suffices to know when  there is no  sequence of glued polyhedra from $$\mathcal{H}(\alpha_1,\ldots,\alpha_i,\alpha_j,\ldots,\alpha_{n+3})$$ to  $$\mathcal{H}(\alpha_1,\ldots,\alpha_j,\alpha_i,\ldots,\alpha_{n+3}).$$ Below are all the possible configurations.
\begin{itemize}
 \item If all the pairs $\alpha_i,\alpha_j$ satisfy $\alpha_i+\alpha_j<\pi$ the polyhedra are glued along  
$\mathcal{H}(\alpha_1,\ldots,\alpha_i+\alpha_j,\ldots,\alpha_{n+3})$;
\item If there is only one pair  $\alpha_i,\alpha_j$ such that $\alpha_i+\alpha_j\geq\pi$, we can glue 
$$\mathcal{H}(\alpha_1,\ldots,\alpha_i,\alpha_j,\alpha_{j+1},\ldots,\alpha_{n+3})$$ with $$\mathcal{H}(\alpha_1,\ldots,\alpha_i,\alpha_{j+1},\alpha_j,\ldots,\alpha_{n+3})$$ and so on, i.e. we can always glue a polyhedron to the one obtained by permuting $\alpha_j$ with the angle at its right. As the list of angles is up to cyclic order,  we arrive at $\mathcal{H}(\alpha_1,\ldots,\alpha_j,\alpha_i,\ldots,\alpha_{n+3})$, hence $\tilde{R}(\alpha_1,\ldots,\alpha_{n+3})$ is connected.
\item  By adapting the argument above with  suitable permutations, it is easy to show that $\tilde{R}(\alpha_1,\ldots,\alpha_{n+3})$ is still connected if there are $\alpha_i,\alpha_j,\alpha_k$ such that   $\alpha_i+\alpha_j\geq\pi$, $\alpha_i+\alpha_k\geq\pi$ and $\alpha_j+\alpha_k<\pi$;
\item In the same way it is easy to show that if $\alpha_i,\alpha_j,\alpha_k$ are as in the statement of the lemma, 
$$\mathcal{H}(\alpha_1,\ldots,\alpha_i,\alpha_j,\alpha_k,\ldots,\alpha_{n+3}) $$
and 
$$\mathcal{H}(\alpha_1,\ldots,\alpha_j,\alpha_i,\alpha_k,\ldots,\alpha_{n+3})  $$
can't be joined and hence $\tilde{R}(\alpha_1,\ldots,\alpha_{n+3})$ is not connected. Moreover every other polyhedron of  $\tilde{R}(\alpha_1,\ldots,\alpha_{n+3})$  can be joined to one of these two polyhedra, and reversing the order of the angles is a bijection between these two components.
\item As the sum of the $\alpha_k$ is equal to $2\pi$ (and $n\geq2$) there can't be a fourth angle $\alpha_l$ such that $\alpha_i,\alpha_j,\alpha_k$ are as in the statement of the lemma and $\alpha_l+\alpha_x\geq\pi$ for $x\in\{1,\ldots,n+3\}$.
\end{itemize}
\end{proof}

We know by the Poincar\'e Theorem \cite[Theorem 4.1]{T} that $R(\alpha_1,\ldots,\alpha_{n+3})$ is isometric to a hyperbolic orbifold if and only if the angle around each singular set of codimension $2$ is  $2\pi/k$ with $k$ integer $>0$. If all these angles are $2\pi$ then the orbifold is a manifold. The singular sets of codimension $2$ in  $R(\alpha_1,\ldots,\alpha_{n+3})$ correspond to codimension $2$ faces of the polyhedra around which facets are glued. There are two possibilities:
\begin{enumerate}
 \item if  $\alpha_k+\alpha_{k+1}<\pi$ and $\alpha_j+\alpha_{j+1}<\pi$ with $2\leq\vert j-k\vert\leq n+1$, then around the codimension 2 face isometric to
$$N:=\mathcal{H}(\alpha_1,\ldots,\alpha_k+\alpha_{k+1},\ldots,\alpha_j+\alpha_{j+1},\ldots,\alpha_{n+3}) $$
 are glued four orthoschemes, corresponding to the four ways of ordering $(\alpha_k,\alpha_{k+1})$ and $(\alpha_j,\alpha_{j+1}).$ As we know that the dihedral angle of each orthoscheme at such codimension $2$ face is $\pi/2$, the total angle around $N$ in $R(\alpha_1,\ldots,\alpha_{n+3})$ is $2\pi$. Hence metrically $N$ is actually not a singular set.
\item if $\alpha_k+\alpha_{k+1}+\alpha_{k+2}<\pi$, then around the codimension 2 face isometric to
$$S:=\mathcal{H}(\alpha_1,\ldots,\alpha_k+\alpha_{k+1}+\alpha_{k+2},\ldots,\alpha_{n+3}) $$
 are glued six orthoschemes corresponding to the six ways of ordering $$(\alpha_k,\alpha_{k+1},\alpha_{k+2}).$$
\end{enumerate}

Let us examine the angle $\theta$ around $S$. It is a sum of  six dihedral angles. Formula~(\ref{eq: cos angle}) gives the cosine of each dihedral angle. It is symmetric in two variables, hence $\theta$ is two times the sum of three different dihedral angles. Moreover
\begin{proposition}\label{lem: calcul angle}
 We have that $\cos(\theta/2)$ is equal to

$$\frac{\sin(\alpha_1)\sin(\alpha_2)\sin(\alpha_3)-\sin(\alpha_1+\alpha_2+\alpha_3)(\sin(\alpha_1)\sin(\alpha_2)+\sin(\alpha_2)\sin(\alpha_3)+\sin(\alpha_3)\sin(\alpha_1))}{\sin(\alpha_1+\alpha_2)\sin(\alpha_2+\alpha_3)\sin(\alpha_3+\alpha_1)}. $$
\end{proposition}

\begin{proof}
This formula is proved in \cite{kojimaal1}, hence we only outline the proof and refer to this reference for more details about the computation. Note that in this reference the result is stated only for $n=2,3$ as they get (\ref{eq: cos angle}) only for these dimensions. The idea is the following. We have to consider the sum of three dihedral angles, and as we know that they are acute, this sum is less that $2\pi$. Hence a gluing of three orthoschemes having those dihedral angles  can be isometrically embedded in the  hyperbolic space (or at least a neighborhood of $S$). The gluing involves four facets glued around  $S$, that gives four outward unit normals $e_1,e_2,e_3,e_4$ spanning a space-like plane in the Minkowski space. Hence we can see these vectors as unit vectors in the Euclidean plane, and the problem is now reduced to find the angle between $e_1$ and $e_4$ knowing the angles between $e_1$ and $e_2$, $e_2$ and $e_3$ and $e_3$ and $e_4$ (the exterior dihedral angles of the orthoschemes).
\end{proof}

Here is a natural question.
\begin{question}
 For which  $(\alpha_1,\ldots,\alpha_{n+3})$ is the cone-manifold $R(\alpha_1,\ldots,\alpha_{n+3})$  isometric to an orbifold?
\end{question}

An intermediate step could be to know if there can exist such orbifolds for all dimensions $n$. This is motivated by the fact that Coxeter orthoschemes and Mostow orbifolds (see Section~\ref{sec: 3}) both don't exist for $n>9$.  
Another analogy is that there exists (at least) $98$ Coxeter orthoschemes (counting one for each infinite family in dimension $\leq 3$) and $94$ Mostow orbifolds (but for dimension $\geq 4$ there are at least $76$ Coxeter orthoschemes and only $28$ Mostow orbifolds). 
But obviously there is no relation between the fact that $R(\alpha_1,\ldots,\alpha_{n+3})$ is an orbifold and the fact that the orthoschemes constituting it are Coxeter. 
For example it is easy to check that $R(\frac{\pi}{2},\frac{\pi}{2},\frac{\pi}{4},\frac{\pi}{4},\frac{\pi}{4},\frac{\pi}{4})$ is a cone-manifold, made of Coxeter orthoschemes (some of them are isometric to the one on the left in Figure~\ref{fig: 2ex cox}).

 Another intermediate step could be to know if there exists non-rational angles $\alpha_i$ such that $R(\alpha_1,\ldots,\alpha_{n+3})$ is an orbifold (in the case of Mostow orbifolds,  the angles have to be rational \cite[3.12]{DM}).

I tried the formula of Proposition \ref{lem: calcul angle} with a computer program, with  $\alpha_i= p\pi/q$, $p$ and $q$ integers, $p<q<100$, as data. A value of the form $\cos(\pi /k)$  (actually $\cos(\pi/2)$)
was reached only for:
\begin{equation}\label{eq: orbifold triple}
 \left(\frac{\pi}{4},\frac{\pi}{4},\frac{5\pi}{12}\right).
\end{equation}
The program may have missed some values and more involved computations would lead too far from the scope of this paper. It is easy to check that if $R(\alpha_1,\ldots,\alpha_{n+3})$ is an orbifold with only (\ref{eq: orbifold triple}) leading to a singular stratum then $n=2$. Two examples of such orbifolds appear in Table~\ref{tableau}, namely $R(\frac{7\pi}{12},\frac{5\pi}{12},\frac{\pi}{2},\frac{\pi}{4},\frac{\pi}{4})$ and $R(\frac{2\pi}{3},\frac{5\pi}{12},\frac{5\pi}{12},\frac{\pi}{4},\frac{\pi}{4}  )$, and using (\ref{eq: cos angle}) we note that the orthoschemes involved in the gluings are not of Coxeter type. 

There is an easy case when 
$R(\alpha_1,\ldots,\alpha_{n+3})$ is a manifold: it is when there does not exist any singular set of codimension $2$, i.e. when the sum of each triple of angles is greater or equal to  $\pi$. It is easy to check that this can happen only for $n=2,3$. If $n=3$  the sum of each triple of angles must be $\pi$, that implies that all the angles are equal. Hence there is only one case, obtained by gluing $60$ times  the orthoscheme shown on the left in Figure~\ref{fig: 2ex cox}. 
For  $n=2$  there are infinitely many  $R(\alpha_1,\ldots,\alpha_{5})$ which are manifolds, for examples the ones obtained by slightly deforming the angles in the list $ \left(\frac{2\pi}{5},\frac{2\pi}{5},\frac{2\pi}{5},\frac{2\pi}{5},\frac{2\pi}{5}\right)$ (they are gluing of right-angled pentagons).

It is proved in \cite{kojimaal1,kojimaal2} that, if $n=2$ or $3$
 and  $R(\alpha_1,\ldots,\alpha_{n+3})$
contains only orthoschemes of type 3, then its
 metric structure   is uniquely determined by the angles $(\alpha_1,\ldots,\alpha_{n+3})$.
 This contrasts with the fact that an orthoscheme can be constructed from  an infinite number of lists of angles. The reason is that permutation of angles and projective transformations do not commute in general.  This leads to the question:
\begin{question}
Does $R(\alpha_1,\ldots,\alpha_{n+3})=R(\alpha_1',\ldots,\alpha_{n+3}')$ mean that $(\alpha_1,\ldots,\alpha_{n+3})=(\alpha_1',\ldots,\alpha_{n+3}')$ (up to the action of the dihedral group)?
\end{question} 

The study of the deformation spaces of polygons with fixed angles appears in \cite{T} as a particular case of the study of the deformation spaces of polyhedra (see next section). It also appears as an exercise in \cite[Problem 2.3.12]{Tlivre}. Detailed studies  can be found in \cite{kojima0,yoshida1,aperyyoshida}  for $n=2$ and in  \cite{aharayamada,MorinHishi} for $n=3$. Both cases are treated in \cite{kojimaal1,kojimaal2}. Note that these references deal mainly  with orthoschemes of type $3$.

Due to the so-called Schwarz--Christoffel map,  $R(\alpha_1,\ldots,\alpha_{n+3})$ can be thought of as a real hyperbolization of the space of configurations of points on the circle, depending on weights $(\alpha_1,\ldots,\alpha_{n+3})$ \cite{T,kojimaal1,IozziPoritz,MorinHishi}. For this reason this construction can be related to many other ones, see e.g. \cite{KM1,IozziPoritz} and references therein. In \cite{KM1} 
it is proved that the spaces of polygons with fixed angles are homeomorphic to the spaces of polygons with fixed edge lengths (so Lemma~\ref{lem: connexite} is the analog of \cite[Theorem 1, Lemma 6]{KM1}). Natural metrics on the moduli spaces of convex polygons in constant curvature Riemannian and Lorentzian spaces are also introduced in \cite{jms}.


\section{Spaces of polyhedra}\label{sec: 3}

\subsection{Configurations of points on the sphere}

In \cite{DM}, the  space  of configurations of points on the sphere with suitable weights is endowed with a complex hyperbolic structure, depending on the weights. It is then possible to find a list of (compact or finite-volume) complex hyperbolic orbifolds. The list was enlarged in \cite{Mostow86} (some of them were known for a long time, we refer to \cite{DM} for more details).  Due to a generalization of the Schwarz--Christoffel map  \cite[8]{T}, the space  of configurations of weighted points on the sphere is homeomorphic to the space of Euclidean cone-metrics on the sphere with prescribed cone-angles (see below) (this also follows from general theorems about the determination of metrics on surfaces by the curvatures \cite{troyanov2,troyanov1}). In \cite{T}, the results of \cite{DM,Mostow86} are recovered by studying such spaces of metrics on the sphere. The two constructions are 
outlined in parallel in  \cite{kojima1}. A bridge between the two constructions is clearly exposed  in   \cite{Troteich}. Moreover this last reference concerns also surfaces of higher genus. We won't review further the numerous works based on \cite{DM} and \cite{T}.

 A \textit{Euclidean metric  with} $N$ \textit{cone singularities of positive curvature}  on the sphere $\mathbb{S}^2$ is a flat metric on $\mathbb{S}^2$ minus $N$ points $x_1,\ldots,x_{N}$, $N\geq 3$, such that a neighborhood of each $x_k$ is isometric to the neighborhood of the apex of a  Euclidean cone of angle $0<\theta_k<2\pi$ (the curvature is $2\pi-\theta_k$).  Such a metric  is uniquely determined up to homotheties by the conformal class of the $N$ punctured sphere and by the numbers $\alpha_k$ (satisfying Gauss--Bonnet condition).
Let $(\alpha_1,\ldots,\alpha_{n+3})$ satisfy (A).
We denote by $C(\alpha_1,\ldots,\alpha_{n+3})$ the set of Euclidean metrics on the sphere with $(n+3)$ cone singularities of positive curvature $2\alpha_i$, up to direct isometries and homotheties. Hence $C(\alpha_1,\ldots,\alpha_{n+3})$ is a connected topological manifold of dimension $2n$. We suppose that the cone-singularities are labeled: if $x_k$ and $x_j$ have the same cone-angle, then exchanging them leads to another metric (our $C(\alpha_1,\ldots,\alpha_{n+3})$ should be written $P(A;2\alpha_1,\ldots,2\alpha_{n+3})$ with the notations of \cite[p. 524]{T}. The notation $C(\alpha_1,\ldots,\alpha_{n+3})$ defined in \cite[p. 524]{T} concerns non-labeled cone-singularities, see remarks after Theorem~\ref{thm: table}). 

Due to the following famous theorem, $C(\alpha_1,\ldots,\alpha_{n+3})$ can also be defined as the space of convex polytopes of $\mathbb{R}^3$ with $n+3$ labeled vertices $x_k$ whose the sum of the angles on the faces around $x_k$  is $2\pi-2\alpha_k$, up to Euclidean direct isometries and homotheties.  In the following we will identify the metric and the polytope.
\begin{theorem}[{Alexandrov Theorem, \cite{alex42,Alexcoll2}}]\label{thm: Alexandrov}
 Let $g$ be a Euclidean metric on the sphere with cone singularities of positive curvature. There exists a convex polytope $P$ in $\mathbb{R}^3$ such that the induced metric on the boundary of $P$ is isometric to $g$. Moreover $P$ is unique up to ambient isometries.
\end{theorem}

It is proved in \cite{KM2} that the spaces of configurations of points on the sphere are also homeomorphic to the spaces of polygons in $\mathbb{R}^3$ with fixed edge lengths. Using a theorem of Minkowski, it is not hard to see that such spaces of polygons are homeomorphic to the spaces of convex polytopes of $\mathbb{R}^3$ with fixed face areas. A duality between the spaces of polygons (in the plane) with fixed angles and the spaces of polygons with fixed edge lengths is proved in \cite{KM1}. In dimension $3$ this duality is expressed between the spaces of polytopes with fixed cone-angles and the spaces of polytopes with fixed face areas. A proof of the   theorem of Minkowski is given in \cite{KlainMinkowski}, simultaneously with a proof of the  Minkowski Inequality, which is the base result to prove the Alexandrov--Fenchel Theorem. A quaternionic structure on the spaces of polygons of $\mathbb{R}^5$ with fixed edge lengths is described in \cite{quaternions}.

To emphasize the analogy between the case of convex polygons in $\mathbb{R}^2$ and the case of convex polytopes in $\mathbb{R}^3$, we describe in the next subsection a local parametrization of $C(\alpha_1,\ldots,\alpha_{n+3})$  close to Thurston's. A difference to his approach is on the choice of unfolding a polytope on the plane. The  use of the
Alexandrov Theorem is  never really necessary, but it simplifies some arguments, and the construction uses elements coming from the original proof by A.~D.~Alexandrov.

\subsection{Polytopes as complex polygons}

A \emph{$N$-gon} of the Euclidean plane is an ordered $N$-tuple of points $(a_1,\ldots,a_N)$ (the vertices) with line segments joining $a_{k-1}$ to $a_k$ (with $a_{N+1}=a_1$). Seeing the Euclidean plane as the complex plane, the set of $N$-gons is identified with $\mathbb{C}^N$.
 Let $P$ be a convex polytope representing an element  of $C(\alpha_1,\ldots,\alpha_{n+3})$. 
We will associate $P$ with   a $(2n+6)$-gon $A(P)$.

\begin{figure}[ht] \begin{center}
\input{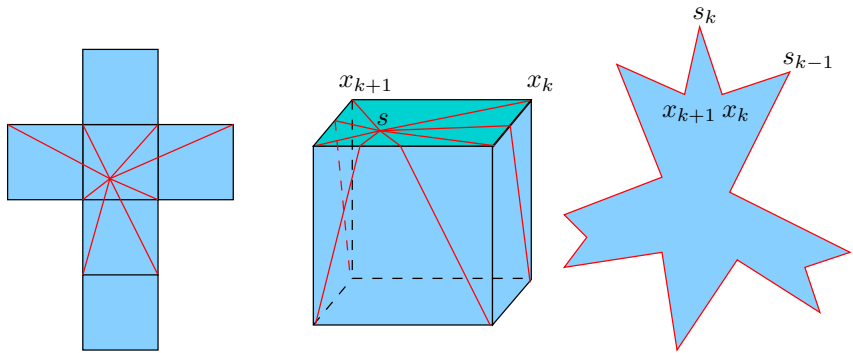}%
\caption{On the right is an Alexandrov unfolding of the cube. On the left is an edge unfolding of the cube, used to determine the shortest geodesics from the source point  to the cone points. \label{fig: Aunfcube}}
\end{center}
\end{figure}

Let us choose a point $s$ on (the boundary of) $P$. The point $s$ is a \emph{source point}. We suppose that $s$ is \emph{generic}: it is not a vertex and for each vertex $x_k$ there is a unique shortest geodesic from $s$ to $x_k$. If we cut $P$ along these shortest geodesics, then it can be unfolded into the plane as a $2(n+3)$-gon: $(n+3)$ vertices are the images of the vertices $x_k$, which alternate with the $(n+3)$ images of $s$. These last ones will be denoted by $(s_1,\ldots,s_{n+3})$. By an abuse of notation the images of the vertices $x_k$ will be still denoted by letters $x_k$, but the vertices of the  $2(n+3)$-gon are labeled such that $s_k$ is between $x_{k}$ and $x_{k+1}$ for the direct order. 
An example  is shown in Figure~\ref{fig: Aunfcube}.

This procedure is known   as \emph{Alexandrov unfolding} or as \emph{star unfolding} (even if the resulting polygon is not necessarily star-shaped, hence we will avoid this last terminology), and  seemingly due to Alexandrov \cite[4.1.2]{AlexCP}, \cite[VI,1]{Alexcoll2}. Actually in those references the source point $s$ is a vertex. With this restriction it is proved that the Alexandrov unfolding is non-overlapping in \cite[Proposition 7.1]{T} (but in Figure~16 of \cite{T} the source point is generic). In our case with $s$ a generic point, the Alexandrov unfolding is also non-overlapping  \cite{AronovORourke}. See also \cite{MP, Pakbook}. Hence $A(P)$ is a simple polygon. 
Alexandrov unfolding is used in \cite{weber} to parameterize the spaces of cone-metrics on the sphere with four cone-singularities of positive curvature and one cone-singularity of negative curvature (i.e.  the angle around the singularity is $>2\pi$). Another way to unfold cone-metrics with five cone-points  in the complex plane is described in \cite{parker}.

By knowing only the $(s_1,\ldots,s_{n+3})$ we can recover $A(P)$ and hence $P$, because the $x_k$ are determined by 
\begin{equation}\label{eq: xk}
 s_{k}-x_k=e^{i2\alpha_{k}}(s_{k-1}-x_k).
\end{equation}

It follows that $A(P)$ is living in a complex vector space of (complex) dimension $(n+3)$. We identify this space with $\mathbb{C}^{n+3}$.
On  $\mathbb{C}^{n+3}$ we define the following Hermitian form:

\begin{equation}\label{eq: M}
M(P,Q)=-\frac{1}{4i}\sum_{k=1}^{n+3}s_k(P)\overline{x_{k}(Q)}-x_{k}(P)\overline{s_k(Q)} + x_{k+1}(P)\overline{s_{k}(Q)}-s_{k}(P)\overline{x_{k+1}(Q)}
\end{equation}
(the signed area of the triangle $0ab$, $(a,b)\in\mathbb{C}^2$, is $\frac{1}{4i}(b\overline{a}-a\overline{b})$).
If $P$ is an element of $C(\alpha_1,\ldots,\alpha_{n+3})$, then $M(P,P)$ is minus the area of $A(P)$ (the face-area of $P$). Here is the analog of Theorem~\ref{thm: alex fenchel}.

\begin{theorem}
  The Hermitian form $M$  has signature $(1,2,n)$.
\end{theorem}

This is proved by an induction on $n$. For $n=0$, polytopes are doubled triangles and the result follows from  Theorem~\ref{thm: alex fenchel} (see the proof of Lemma~\ref{lem: isom embed}). This can also be seen directly by the facts that the triangles have positive area and that the area is invariant under translations.

For any $n$, we go back to $n-1$ by the process called ``cutting and gluing'' and seemingly due to Alexandrov \cite[Lemma 1,p 226]{Alexcoll2}, see also \cite[17.5]{busemann} and \cite[Proposition 3.3]{T}. Cutting and gluing is as follows: if two cone points, with curvatures $2\alpha_k$ and $2\alpha_j$ such that $\alpha_k+\alpha_j<\pi$, are sufficiently close (such points always exist), then we  cut the geodesic joining them. To the two resulting geodesics it is possible to glue a Euclidean cone of curvature $2(\alpha_k+\alpha_j)$ in such a way that the singularities at $x_k$ and $x_j$ disappear. The area of the old metric is the area of the new metric minus the area of the cone.  The similar procedure applied to polygons instead of polytopes  is the way used in \cite{BG} to prove Theorem~\ref{thm: alex fenchel} of the current paper.

The Hermitian form (\ref{eq: M}) can be considered as the complex mixed-area. It is a classical form on the space of the $N$-gons, see e.g. \cite{perpendicularpolygons} and the references therein. In  \cite{perpendicularpolygons} the form  is related to a larger family of Hermitian forms. Moreover polygons are considered as finite Fourier Series. It is possible to prove the Alexandrov--Fenchel Theorem (actually the Minkowski Inequality) for convex curves using Fourier Series. Perhaps it is possible to compute the signature of $M$ using this point of view. For more details we refer to  
\cite{groemer}, especially Formula (4.3.3) and Remarks and References of 4.1. 

The quotient of  $\mathbb{C}^{n+3}$ by the kernel of $M$ describes the unfoldings up to  translations. 
This quotient is a complex vector space of (complex) dimension $(n+1)$.
An Alexandrov unfolding is clearly a well-defined and injective map from a (sufficiently small) neighborhood $U$ of $P$ in $C(\alpha_1,\ldots,\alpha_{n+3})$ if the unfoldings are moreover considered up to rotations and homotheties. Hence $U$ is mapped homeomorphically to a subset of the set of negative vectors (for a Hermitian form of signature $(1,n)$) of the quotient of a complex vector space of dimension $(n+1)$  by complex conjugation:   Alexandrov unfolding provides charts from $C(\alpha_1,\ldots,\alpha_{n+3})$ to $\mathbb{CH}^n$, the complex hyperbolic space of (complex) dimension $n$  --- we refer to \cite{epstein,Goldman} for details about $\mathbb{CH}^n$.

Let $P$ be an element of $C(\alpha_1,\ldots,\alpha_{n+3})$. Let $T$ be a geodesic triangulation of $P$ 
such that the cone-points are exactly the vertices of $T$ (we will call such a triangulation a \emph{cone-triangulation}). The fact that cone-triangulations exist is obvious if we use the Alexandrov Theorem (it suffices to triangulate the faces of the corresponding convex polytope), but precisely the existence of cone-triangulations  is a step in the proof of this theorem, see  \cite[p. 130]{busemann}, \cite[Proposition 3.1]{T}. If we cut along some edges of $T$ it is possible to unfold $P$ to the complex plane: this is an \emph{edge unfolding} (it is not necessarily non-overlapping). To each edge of $T$ is associated a complex number, and $n+1$ complex numbers suffice to recover $P$ \cite[Proposition 3.2]{T}: edge unfolding provides another charts from $C(\alpha_1,\ldots,\alpha_{n+3})$ to $\mathbb{CH}^n$. We check that these two kinds of local coordinates are compatible.

\begin{lemma}\label{lem: AU2tri}
 Let $P$ be in $C(\alpha_1,\ldots,\alpha_{n+3})$ and   let $A$ be an Alexandrov unfolding of $P$.
Let $U$ be a neighborhood of $P$. If $U$ is  sufficiently small, then there exists an edge unfolding  $E$ of $P$ such that $A$ and $E$ are homeomorphisms on $U$ and such that  there is an isometric linear bijection sending $A(U)$ to $E(U)$.
\end{lemma}

\begin{figure}[ht] \begin{center}
\input{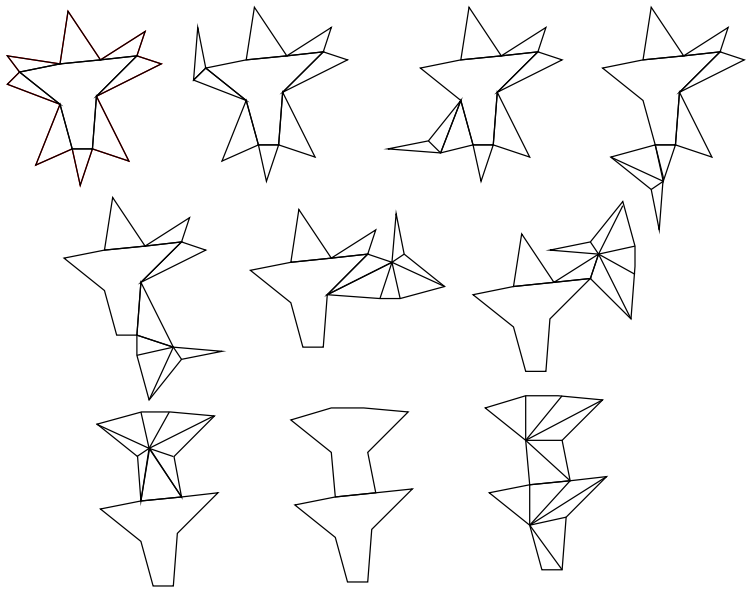}%
\caption{From an Alexandrov unfolding  to an edge unfolding of the cube (example for the proof of Lemma~\ref{lem: AU2tri}).\label{fig: AU2triangles}}
\end{center}
\end{figure}
 \begin{proof}
  The Alexandrov unfolding  $A(P)$ contains a simple polygon $SA(P)$ with vertices $x_1,\ldots,x_{n+3}$ (it is simple because $A(P)$ is). To go from $SA(P)$ to $A(P)$ one has to add to each edge $x_{k}x_{k+1}$ of $SA(P)$ the triangle $T_k:=x_{k}s_kx_{k+1}$. We will ``roll'' all the triangles $T_k$ around $SA(P)$. More precisely, we perform on $T_2$ a rotation of angle $2\alpha_3$ and center $x_3$. $T_2$ is now glued on the edge $x_3s_3$. We rotate the union of $T_2$ and $T_3$   around $x_4$ by an angle $2\alpha_4$, and so on. At the end all the $T_k$ are glued around $T_1$ (all the $s_k$ go to $s_1$), see Figure~\ref{fig: AU2triangles} for an example with the cube. The gluing of all the $T_k$ around $s_k$ gives a simple polygon (because the sum of the angles around $s$ is $2\pi$). So we get  two simple polygons glued along the edge $x_1x_2$. A triangulation of each of them is exactly an edge unfolding of  $P$ (note that there is no reason for the union to be simple). Let us denote by $Q$ this union.

The coordinates of $E(U)$ are  vectors associated to the diagonals of $Q$, namely the differences of the coordinates of their endpoints. These endpoints are the vertices of $Q$. By the preceding paragraph, these vertices are obtained from $x_1,\ldots,x_{n+3}$ by compositions of rotations, and by (\ref{eq: xk}), the $x_i$ are linear functions of the $s_i$, which are the coordinates of $A(U)$. This describes a linear map sending $A(U)$ to $E(U)$.
Obviously this operation  preserves the area and is bijective.
 \end{proof}

From now we summarize the results of \cite{T}. For the coordinates given by the  edge unfoldings on $C(\alpha_1,\ldots,\alpha_{n+3})$, the changes of charts correspond to flipping the edges of $T$ (when two triangles form a quadrilateral with a diagonal, to \textit{flip} is to delete the diagonal and to choose the other). In terms of the complex coordinates, the changes of charts are  linear maps, and obviously  isometries. 
This gives a structure of complex hyperbolic manifold of complex dimension $n$ on $C(\alpha_1,\ldots,\alpha_{n+3})$.
This manifold   is not complete (as cone points can  collide). We  denote by $\overline{C}(\alpha_1,\ldots,\alpha_{n+3})$  its metric completion. Then $\overline{C}(\alpha_1,\ldots,\alpha_{n+3})$ has a structure of complex hyperbolic cone manifold (a cone manifold structure for non constant curvature  is less obvious to define than in the constant curvature case. We refer to \cite{T} for a precise definition). 

The collusion of  two cone-points  $x_k$ and $x_j$ describes  a singular stratum of (complex) codimension $1$ in $\overline{C}(\alpha_1,\ldots,\alpha_{n+3})$ (the collusion is possible only if $\alpha_k+\alpha_j<\pi$). The main point is that the singular curvature around the stratum  is $2\alpha_k+2\alpha_j$ \cite[Proposition 3.5]{T}. Hence it is easy to know for which $(\alpha_1,\ldots,\alpha_{n+3})$ the cone-angles around the (real) codimension $2$ strata are of the form $2\pi / k$. By the  Poincar\'e Theorem \cite[Theorem 4.1]{T} this means that $\overline{C}(\alpha_1,\ldots,\alpha_{n+3})$ is a complex hyperbolic orbifold. There exists $36$ of these orbifolds, which are listed in Table~\ref{tableau}.

Moreover $\overline{C}(\alpha_1,\ldots,\alpha_{n+3})$  has finite volume and it is compact if and only if there is no subset of $(\alpha_1,\ldots,\alpha_{n+3})$ summing to $\pi$ \cite[Proof of Theorem~0.2]{T}.

\section{Spaces of polygons into spaces of polyhedra}\label{sec: 4}

The set of convex polytopes contains  \emph{degenerated (convex) polytopes} which are obtained by ``doubling'' a convex polygon. Doubling is  gluing isometrically along the edges a polygon to its image by a reflection in a line. Such a reflection reverses the labeling of the angles, and hence 
 there is a canonical injection $f$ from  $R(\alpha_1,\ldots,\alpha_{n+3})$ to $\overline{C}(\alpha_1,\ldots,\alpha_{n+3})$. The metric structure on each space is given by the face-area, so it is not surprising that $f$ is an isometry.

\begin{lemma}\label{lem: isom embed}
For each choice of an order on $(\alpha_1,\ldots,\alpha_{n+3})$, the composition of $f$ with an Alexandrov unfolding   extends to an isometric linear map from $(\mathbb{R}^{n+1},2m)$ to $(\mathbb{C}^{n+1},M)$. 
\end{lemma}
\begin{proof}
 Let $P$ be a convex polygon with exterior angles $(\alpha_1,\ldots,\alpha_{n+3})$ and let  $(u_1,\ldots,u_{n+3})$ be its set of outward unit normals. We choose a point $s$ in the interior of $P$. Without loss of generality, let us suppose that $s$ is the origin.
We denote by $s_k$ the image of $s$ by a reflection in the edge $x_{k}x_{k+1}$. The polygon $x_1s_1x_2\ldots x_{n+3}s_{n+3}$ is an Alexandrov unfolding of the doubling of $P$, and then the  linear extension of $f$ is the following map from $\mathbb{R}^{n+1}$ to $\mathbb{C}^{n+1}$:
$$(h_1,\ldots,h_{n+3})\mapsto (2h_1u_1,\ldots,2h_{n+3}u_{n+3}).$$
It follows that $M(f(u_k),f(u_j))=0$ if $2\leq \vert j-k\vert\leq n+1$, and writing (\ref{eq: xk}) as
$$x_k=\frac{1}{2i}\frac{1}{\sin(\alpha_{k})}\left( e^{i\alpha_{k}}s_{k-1}-e^{-i\alpha_{k}}s_{k}\right) $$
we compute easily that $M$ is two times  $m$ (compare with (\ref{eq:squarrednorm}) and (\ref{eq: mui})):  
\begin{eqnarray*}
 \ M(f(u_k),f(u_{k+1}))&=& 4M(u_k,u_{k+1}) \\
\ &=& -\frac{1}{i}\left(-u_k\overline{x_{k+1}(u_{k+1})}-x_{k+1}(u_k)\overline{u_{k+1}}\right) \\
\ &=& \frac{1}{i} \left(\frac{1}{2i}\frac{1}{\sin(\alpha_{k+1})}+\frac{1}{2i}\frac{1}{\sin(\alpha_{k+1})} \right) \\
\ &=& -\frac{1}{\sin(\alpha_{k+1})}=2m(u_k,u_{k+1}),
\end{eqnarray*}
\begin{eqnarray*}
 \ M(f(u_k),f(u_k))&=& 4 M(u_k,u_k)\\
\ &=& -\frac{1}{i}\left(u_k\overline{x_{k}(u_k)}-
\overline{u_k}x_{k}(u_k)+\overline{u_k}x_{k+1}(u_k)- u_k\overline{x_{k+1}(u_k)} \right) \\
\ &=&-\frac{1}{i}\left(\frac{1}{i}\frac{\cos(\alpha_k)}{\sin(\alpha_k)}+\frac{1}{i}\frac{\cos(\alpha_{k+1})}{\sin(\alpha_{k+1})}\right) \\
\ &=& \frac{\sin(\alpha_{k}+\alpha_{k+1})}{\sin(\alpha_{k})\sin(\alpha_{k+1})}=2m(u_k,u_k). 
\end{eqnarray*}
\end{proof}

It follows that on the image of $\mathbb{R}^{n+1}$
 in $\mathbb{C}^{n+1}$, the Hermitian form $M$ has real values. Moreover this image has maximal real dimension, hence it is a real form of $\mathbb{C}^{n+1}$.
To each real form is associated a unique real structure (= anti-linear involution) compatible with the Hermitian structure, whose fixed-points set is the real form. Here the real structure corresponds exactly to the complex conjugation. We follow \cite{Goldman} for the definitions and refer to it for more details. 

The real structure on charts given by polygons comes from a global isometric involution reversing the orientation  on $\overline{C}(\alpha_1,\ldots,\alpha_{n+3})$, denoted by $\rho$.  The involution $\rho$ can be  described as the reflection of polytopes in a plane. As the vertices  are labeled, the fixed-points set of $\rho$ is exactly  
$R(\alpha_1,\ldots,\alpha_{n+3})$. We get easily the metric structure of this set for the orbifolds discovered by Deligne and Mostow.

\begin{theorem}\label{thm: table}
 Table~\ref{tableau} gives the metric structure of $R(\alpha_1,\ldots,\alpha_{n+3})$ for Deligne--Mostow orbifolds.
\end{theorem}

\begin{proof}
If $R(\alpha_1,\ldots,\alpha_{n+3})$ has no singular set of codimension $2$, then it  is a manifold. This occurs when for each triple $(\alpha_i,\alpha_j,\alpha_k)$,  $\alpha_i+\alpha_j+\alpha_k\geq\pi$. 
The two cases marked as orbifolds in Table~\ref{tableau} have only one singular stratum, which is represented by the triple $(\frac{\pi}{4},\frac{\pi}{4},\frac{5\pi}{12})$ (see (\ref{eq: orbifold triple})). 

From the discussion below Proposition~\ref{lem: calcul angle}, we know that the other examples are neither manifolds nor orbifolds. We check this fact.  If there exists three angles $\alpha_i,\alpha_j,\alpha_k$ such that the dihedral angles given by  $(\alpha_i,\alpha_j,\alpha_k)$ and $(\alpha_j,\alpha_i,\alpha_k)$ are $\pi / 4$, then the total angle around the singular set defined by these angles
  is at least four times $\pi/4$, plus something less than $\pi$ (we know that these dihedral angles of orthoschemes are $<\pi/2$), hence it can't be  $2\pi /k$. It is easy to see that such a condition occurs  for the triples $ (\frac{\pi}{4},\frac{\pi}{4},\frac{\pi}{4})$, $(\frac{\pi}{3},\frac{\pi}{3},\frac{\pi}{6})$, $(\frac{3\pi}{8},\frac{3\pi}{8},\frac{\pi}{8})$, $(\frac{\pi}{10},\frac{2\pi}{5},\frac{2\pi}{5})$, $(\frac{\pi}{12},\frac{5\pi}{12},\frac{5\pi}{12})$ and $(\frac{\pi}{18},\frac{4\pi}{9},\frac{4\pi}{9})$, which cover all the  cases indicated in Table~\ref{tableau}.

\end{proof}

\renewcommand{\arraystretch}{2}
\begin{center}\begin{table}
\begin{tabular}{|c|c|c||c|c|c|}
\hline  \textbf{T}&\textbf{Angles}&\textbf{S }  & \textbf{T}&\textbf{Angles}&\textbf{S }\\ 
\hline&Dimension $5$&& 43  & $\displaystyle\frac{\pi}{2},\frac{\pi}{2},\frac{\pi}{3},\frac{\pi}{3},\frac{\pi}{3}$&M\\
 3 & $\displaystyle\frac{\pi}{4},\frac{\pi}{4},\frac{\pi}{4},\frac{\pi}{4},\frac{\pi}{4},\frac{\pi}{4},\frac{\pi}{4},\frac{\pi}{4}$&C 
& 45  & $\displaystyle\displaystyle\frac{3\pi}{4},\frac{\pi}{8},\frac{3\pi}{8},\frac{3\pi}{8},\frac{3\pi}{8}$&C \\
&Dimension $4$&& 46  & $\displaystyle\frac{5\pi}{8},\frac{5\pi}{8},\frac{\pi}{4},\frac{\pi}{4},\frac{\pi}{4}$&C \\
 4& $\displaystyle\frac{\pi}{2},\frac{\pi}{4},\frac{\pi}{4},\frac{\pi}{4},\frac{\pi}{4},\frac{\pi}{4},\frac{\pi}{4}$&C 
& 47  & $\displaystyle\frac{\pi}{2},\frac{3\pi}{8},\frac{3\pi}{8},\frac{3\pi}{8},\frac{3\pi}{8}$&M\\
&Dimension $3$&&48 & $\displaystyle\frac{2\pi}{9},\frac{4\pi}{9},\frac{4\pi}{9},\frac{4\pi}{9},\frac{4\pi}{9}$&M\\
1 & $\displaystyle\frac{\pi}{3},\frac{\pi}{3},\frac{\pi}{3},\frac{\pi}{3},\frac{\pi}{3},\frac{\pi}{3}$&M
& 49  & $\displaystyle\frac{\pi}{10},\frac{7\pi}{10},\frac{2\pi}{5},\frac{2\pi}{5},\frac{2\pi}{5}$&C \\
 5 & $\displaystyle\frac{3\pi}{4},\frac{\pi}{4},\frac{\pi}{4},\frac{\pi}{4},\frac{\pi}{4},\frac{\pi}{4}$&C 
& 57  & $\displaystyle\frac{2\pi}{3},\frac{\pi}{12},\frac{5\pi}{12},\frac{5\pi}{12},\frac{5\pi}{12}$&C \\
6 & $\displaystyle\frac{\pi}{2},\frac{\pi}{2},\frac{\pi}{4},\frac{\pi}{4},\frac{\pi}{4},\frac{\pi}{4}$&C 
& 65 & $\displaystyle\frac{7\pi}{12},\frac{7\pi}{12},\frac{\pi}{6},\frac{\pi}{3},\frac{\pi}{3}$&C \\
 39 & $\displaystyle\frac{\pi}{2},\frac{\pi}{6},\frac{\pi}{3},\frac{\pi}{3},\frac{\pi}{3},\frac{\pi}{3}$&C 
& 68  & $\displaystyle\frac{5\pi}{6},\frac{5\pi}{12},\frac{\pi}{4},\frac{\pi}{4},\frac{\pi}{4}$&C\\
44& $\displaystyle\displaystyle\frac{\pi}{8},\frac{3\pi}{8},\frac{3\pi}{8},\frac{3\pi}{8},\frac{3\pi}{8},\frac{3\pi}{8}$&C 
& 69 & $\displaystyle\frac{2\pi}{3},\frac{7\pi}{12},\frac{\pi}{4},\frac{\pi}{4},\frac{\pi}{4}$&C\\
66 & $\displaystyle\frac{7\pi}{12},\frac{5\pi}{12},\frac{\pi}{4},\frac{\pi}{4},\frac{\pi}{4},\frac{\pi}{4}$&C 
& 70 & $\displaystyle\frac{2\pi}{3},\frac{5\pi}{12},\frac{5\pi}{12},\frac{\pi}{4},\frac{\pi}{4}$&O\\
67& $\displaystyle\frac{5\pi}{12},\frac{5\pi}{12},\frac{5\pi}{12},\frac{\pi}{4},\frac{\pi}{4},\frac{\pi}{4}$&C 
& 71 & $\displaystyle\frac{7\pi}{12},\frac{5\pi}{12},\frac{\pi}{2},\frac{\pi}{4},\frac{\pi}{4}$&O\\
&Dimension $2$&& 72 & $\displaystyle\frac{\pi}{2},\frac{\pi}{4},\frac{5\pi}{12},\frac{5\pi}{12},\frac{5\pi}{12}$&M\\
2 & $\displaystyle\frac{2\pi}{3},\frac{\pi}{3},\frac{\pi}{3},\frac{\pi}{3},\frac{\pi}{3}$&M
& 73 & $\displaystyle\frac{7\pi}{12},\frac{5\pi}{12},\frac{\pi}{3},\frac{\pi}{3},\frac{\pi}{3}$&M \\
7 &$\displaystyle\frac{\pi}{2},\frac{3\pi}{4},\frac{\pi}{4},\frac{\pi}{4},\frac{\pi}{4}$&C
& 74 & $\displaystyle\frac{\pi}{2},\frac{\pi}{3},\frac{\pi}{3},\frac{5\pi}{12},\frac{5\pi}{12}$&M \\
8 & $\displaystyle\frac{\pi}{2},\frac{\pi}{2},\frac{\pi}{2},\frac{\pi}{4},\frac{\pi}{4}$&M 
& 75 & $\displaystyle\frac{\pi}{3},\frac{5\pi}{12},\frac{5\pi}{12},\frac{5\pi}{12},\frac{5\pi}{12}$&M\\
9&  $\displaystyle\frac{2\pi}{5},\frac{2\pi}{5},\frac{2\pi}{5},\frac{2\pi}{5},\frac{2\pi}{5}$&M
& 78  & $\displaystyle\frac{4\pi}{15},\frac{8\pi}{15},\frac{2\pi}{5},\frac{2\pi}{5},\frac{2\pi}{5}$&M \\
40 & $\displaystyle\frac{5\pi}{6},\frac{\pi}{6},\frac{\pi}{3},\frac{\pi}{3},\frac{\pi}{3}$&C
& 79& $\displaystyle\frac{\pi}{18},\frac{11\pi}{18},\frac{4\pi}{9},\frac{4\pi}{9},\frac{4\pi}{9}$&C \\
41&  $\displaystyle\frac{2\pi}{3},\frac{\pi}{3},\frac{\pi}{3},\frac{\pi}{2},\frac{\pi}{6}$&C 
& 85 &  $\displaystyle\frac{7\pi}{10},\frac{11\pi}{20},\frac{\pi}{4},\frac{\pi}{4},\frac{\pi}{4}$&C\\
42&  $\displaystyle\frac{\pi}{2},\frac{\pi}{2},\frac{\pi}{2},\frac{\pi}{3},\frac{\pi}{6}$&M 
& 89 &  $\displaystyle\frac{7\pi}{12},\frac{7\pi}{24},\frac{3\pi}{8},\frac{3\pi}{8},\frac{3\pi}{8}$&M\\
\hline
\end{tabular}
 \caption{\label{tableau} The angles $(\alpha_1,\ldots,\alpha_{n+3})$ are those for which $\overline{C}(\alpha_1,\ldots,\alpha_{n+3})$ is a complex hyperbolic orbifold, given by the list in \cite{DM}. The column named \textbf{T} gives the number of the orbifold in the list of \cite{T}. The column named \textbf{S} gives the structure of $R(\alpha_1,\ldots,\alpha_{n+3})$: M means that it is a manifold, O that it is an orbifold and C that it is a cone-manifold.}
\end{table}\end{center}

It is also possible to study the spaces of cone-metrics without labeling  the cone-points. In this case there are more orbifolds: if $\alpha_k=\alpha_j$ and $\alpha_k+\alpha_j<\pi$, then the angle around the stratum is half of the one obtained with labeling.  Such orbifolds were founded in \cite{Mostow86,Mostow88}. The  complete list is achieved in \cite{T}
(the list is known to be complete due to \cite[3.12]{DM} and \cite{felikson}). This list  contains the list founded in \cite{DM} and given in Table~\ref{tableau}.
In this case of non-labeling of the cone-points, the fixed-points set of $\rho$ contains the spaces of polygons and some polytopes obtained by doubling convex caps (which can be seen as convex isometric embeddings of Euclidean metrics with conical singularities on the closed disc). Answering the following question should be a step in the study of the fixed-points set of $\rho$ for Mostow orbifolds.

\begin{question}
 Is it possible to describe a (real) hyperbolic structure on the space of convex caps with fixed cone-angles?
\end{question}

The following works concern real forms of  complex hyperbolic orbifolds, with approaches different from ours:
 \cite{aperyyoshida} for $n=2$,  \cite{ACT3,ACT4} for $n=3$, 
 \cite{yoshida2} and \cite{ACT2} (announced in \cite{ACT1}) for $n=4$,  \cite{chu-2007} for $n=5$, and the references therein. From the fixed-points sets appearing in those works the following question arises:

\begin{question}
 Is it possible to describe a (real) hyperbolic structure on the space of centrally symmetric convex polytopes with fixed cone-angles?
\end{question}

\bibliographystyle{alpha}

\begin{thebibliography}{ACT07b}

\bibitem[ACT03]{ACT1}
D.~Allcock, J.~A. Carlson, and D.~Toledo.
\newblock Real cubic surfaces and real hyperbolic geometry.
\newblock {\em C. R. Math. Acad. Sci. Paris}, 337(3):185--188, 2003.

\bibitem[ACT06]{ACT4}
D.~Allcock, J.~A. Carlson, and D.~Toledo.
\newblock Nonarithmetic uniformization of some real moduli spaces.
\newblock {\em Geom. Dedicata}, 122:159--169, 2006.

\bibitem[ACT07a]{ACT2}
D.~Allcock, J.~A. Carlson, and D.~Toledo.
\newblock Hyperbolic geometry and moduli of real cubic surfaces.
\newblock To appear in \textit{Ann. Ecole Norm. Sup.}, 2007.

\bibitem[ACT07b]{ACT3}
D.~Allcock, J.~A. Carlson, and D.~Toledo.
\newblock Hyperbolic geometry and the moduli space of real binary sextics.
\newblock In {\em Arithmetic and geometry around hypergeometric functions},
  volume 260 of {\em Progr. Math.}, pages 1--22. Birkh\"auser, Basel, 2007.

\bibitem[Ale37]{AlexAF}
A.~D. Alexandrov.
\newblock On the theory of mixed volumes {II}.
\newblock {\em Mat. Sbornik}, 44:1205--1238, 1937.
\newblock (Russian. Translated in \cite{Alexcoll1}).

\bibitem[Ale42]{alex42}
A.~D. Alexandroff.
\newblock Existence of a convex polyhedron and of a convex surface with a given
  metric.
\newblock {\em Rec. Math. [Mat. Sbornik] N.S.}, 11(53):15--65, 1942.
\newblock (Russian. Translated in \cite{Alexcoll1}).

\bibitem[Ale96]{Alexcoll1}
A.~D. Alexandrov.
\newblock {\em Selected works. {P}art {I}}, volume~4 of {\em Classics of Soviet
  Mathematics}.
\newblock Gordon and Breach Publishers, Amsterdam, 1996.
\newblock Selected scientific papers, Translated from the Russian by P. S. V.
  Naidu, Edited and with a preface by Yu. G. Reshetnyak and S. S. Kutateladze.

\bibitem[Ale05]{AlexCP}
A.~D. Alexandrov.
\newblock {\em Convex polyhedra}.
\newblock Springer Monographs in Mathematics. Springer-Verlag, Berlin, 2005.
\newblock Translated from the 1950 Russian edition by N. S. Dairbekov, S. S.
  Kutateladze and A. B. Sossinsky, With comments and bibliography by V. A.
  Zalgaller and appendices by L. A. Shor and Yu. A. Volkov.

\bibitem[Ale06]{Alexcoll2}
A.~D. Alexandrov.
\newblock {\em A. {D}. {A}lexandrov selected works. {P}art {II}}.
\newblock Chapman \& Hall/CRC, Boca Raton, FL, 2006.
\newblock Intrinsic geometry of convex surfaces, Edited by S. S. Kutateladze,
  Translated from the Russian by S. Vakhrameyev.

\bibitem[AO92]{AronovORourke}
B.~Aronov and J.~O'Rourke.
\newblock Nonoverlap of the star unfolding.
\newblock {\em Discrete Comput. Geom.}, 8(3):219--250, 1992.
\newblock ACM Symposium on Computational Geometry (North Conway, NH, 1991).

\bibitem[AY98]{aperyyoshida}
F.~Ap{\'e}ry and M.~Yoshida.
\newblock Pentagonal structure of the configuration space of five points in the
  real projective line.
\newblock {\em Kyushu J. Math.}, 52(1):1--14, 1998.

\bibitem[AY99]{aharayamada}
K.~Ahara and K.~Yamada.
\newblock Shapes of hexagrams.
\newblock {\em J. Math. Sci. Univ. Tokyo}, 6(3):539--558, 1999.

\bibitem[BG92]{BG}
C.~Bavard and {\'E}.~Ghys.
\newblock Polygones du plan et poly\`edres hyperboliques.
\newblock {\em Geom. Dedicata}, 43(2):207--224, 1992.

\bibitem[BI08]{BobenkoIzmestiev}
A.~I. Bobenko and I.~Izmestiev.
\newblock Alexandrov's theorem, weighted {D}elaunay triangulations, and mixed
  volumes.
\newblock {\em Ann. Inst. Fourier (Grenoble)}, 58(2):447--505, 2008.

\bibitem[Bus58]{busemann}
H.~Busemann.
\newblock {\em Convex surfaces}.
\newblock Interscience Tracts in Pure and Applied Mathematics, no. 6.
  Interscience Publishers, Inc., New York, 1958.

\bibitem[Chu07]{chu-2007}
K.~Chu.
\newblock On the geometry of the moduli space of real binary octics, 2007.
\newblock arXiv.org:0708.0419.

\bibitem[Deb90]{DB}
H.~E. Debrunner.
\newblock Dissecting orthoschemes into orthoschemes.
\newblock {\em Geom. Dedicata}, 33(2):123--152, 1990.

\bibitem[DM86]{DM}
P.~Deligne and G.~D. Mostow.
\newblock Monodromy of hypergeometric functions and nonlattice integral
  monodromy.
\newblock {\em Inst. Hautes \'Etudes Sci. Publ. Math.}, (63):5--89, 1986.

\bibitem[Eps87]{epstein}
D.~B.~A. Epstein.
\newblock Complex hyperbolic geometry.
\newblock In {\em Analytical and geometric aspects of hyperbolic space
  ({C}oventry/{D}urham, 1984)}, volume 111 of {\em London Math. Soc. Lecture
  Note Ser.}, pages 93--111. Cambridge Univ. Press, Cambridge, 1987.

\bibitem[Fel97]{felikson}
A.~A. Felikson.
\newblock On {T}hurston signatures.
\newblock {\em Uspekhi Mat. Nauk}, 52(4(316)):217--218, 1997.
\newblock Translation in Russian Math. Surveys 52 (1997), no. 4, 826--827.

\bibitem[FI]{AF}
F.~Fillastre and I.~Izmestiev.
\newblock Construct hyperbolic polyhedra from {A}lexandrov--{F}enchel theorem.
\newblock Provisional title. Work in progress.

\bibitem[FL04]{quaternions}
P.~Foth and G.~Lozano.
\newblock The geometry of polygons in {$\Bbb R\sp 5$} and quaternions.
\newblock {\em Geom. Dedicata}, 105:209--229, 2004.

\bibitem[FRS85]{perpendicularpolygons}
J.~Chris Fisher, D.~Ruoff, and J.~Shilleto.
\newblock Perpendicular polygons.
\newblock {\em Amer. Math. Monthly}, 92(1):23--37, 1985.

\bibitem[Gol99]{Goldman}
W.~M. Goldman.
\newblock {\em Complex hyperbolic geometry}.
\newblock Oxford Mathematical Monographs. The Clarendon Press Oxford University
  Press, New York, 1999.
\newblock Oxford Science Publications.

\bibitem[Gro96]{groemer}
H.~Groemer.
\newblock {\em Geometric applications of {F}ourier series and spherical
  harmonics}, volume~61 of {\em Encyclopedia of Mathematics and its
  Applications}.
\newblock Cambridge University Press, Cambridge, 1996.

\bibitem[IH85]{ImHof2}
H.-C. Im~Hof.
\newblock A class of hyperbolic {C}oxeter groups.
\newblock {\em Exposition. Math.}, 3(2):179--186, 1985.

\bibitem[IH90]{ImHof1}
H.-C. Im~Hof.
\newblock Napier cycles and hyperbolic {C}oxeter groups.
\newblock {\em Bull. Soc. Math. Belg. S\'er. A}, 42(3):523--545, 1990.
\newblock Algebra, groups and geometry.

\bibitem[IP99]{IozziPoritz}
A.~Iozzi and J.~A. Poritz.
\newblock The moduli space of boundary compactifications of {${\rm SL}(2,{\bf
  R})$}.
\newblock {\em Geom. Dedicata}, 76(1):65--79, 1999.

\bibitem[Izm08]{Izmestiev}
I.~Izmestiev.
\newblock The colin de {V}erdi\`ere number and graphs of polytopes.
\newblock To appear {\em Israel J. Math.} arXiv:0704.0349, 2008.

\bibitem[Kla04]{KlainMinkowski}
D.~A. Klain.
\newblock The {M}inkowski problem for polytopes.
\newblock {\em Adv. Math.}, 185(2):270--288, 2004.

\bibitem[KM95]{KM1}
M.~Kapovich and J.~Millson.
\newblock On the moduli space of polygons in the {E}uclidean plane.
\newblock {\em J. Differential Geom.}, 42(2):430--464, 1995.

\bibitem[KM96]{KM2}
M.~Kapovich and J.~Millson.
\newblock The symplectic geometry of polygons in {E}uclidean space.
\newblock {\em J. Differential Geom.}, 44(3):479--513, 1996.

\bibitem[KNY99]{kojimaal1}
S.~Kojima, H.~Nishi, and Y.~Yamashita.
\newblock Configuration spaces of points on the circle and hyperbolic {D}ehn
  fillings.
\newblock {\em Topology}, 38(3):497--516, 1999.

\bibitem[Koj01]{kojima1}
S.~Kojima.
\newblock Complex hyperbolic cone structures on the configuration spaces.
\newblock {\em Rend. Istit. Mat. Univ. Trieste}, 32(suppl. 1):149--163 (2002),
  2001.
\newblock Dedicated to the memory of Marco Reni.

\bibitem[KY93]{kojima0}
S.~Kojima and Y.~Yamashita.
\newblock Shapes of stars.
\newblock {\em Proc. Amer. Math. Soc.}, 117(3):845--851, 1993.

\bibitem[MN00]{MorinHishi}
B.~Morin and H.~Nishi.
\newblock Hyperbolic structures on the configuration space of six points in the
  projective line.
\newblock {\em Adv. Math.}, 150(2):202--232, 2000.

\bibitem[Mos86]{Mostow86}
G.~D. Mostow.
\newblock Generalized {P}icard lattices arising from half-integral conditions.
\newblock {\em Inst. Hautes \'Etudes Sci. Publ. Math.}, (63):91--106, 1986.

\bibitem[Mos88]{Mostow88}
G.~D. Mostow.
\newblock On discontinuous action of monodromy groups on the complex
  {$n$}-ball.
\newblock {\em J. Amer. Math. Soc.}, 1(3):555--586, 1988.

\bibitem[MP08]{MP}
E.~Miller and I.~Pak.
\newblock Metric combinatorics of convex polyhedra: cut loci and nonoverlapping
  unfoldings.
\newblock {\em Discrete Comput. Geom.}, 39(1-3):339--388, 2008.

\bibitem[Pak08]{Pakbook}
I.~Pak.
\newblock Lectures on discrete and polyhedral geometry.
\newblock To be published. Preliminary version available at author's web page,
  2008.

\bibitem[Par06]{parker}
J.~R. Parker.
\newblock Cone metrics on the sphere and {L}ivn\'e's lattices.
\newblock {\em Acta Math.}, 196(1):1--64, 2006.

\bibitem[Sch93]{schneider}
R.~Schneider.
\newblock {\em Convex bodies: the {B}runn-{M}inkowski theory}, volume~44 of
  {\em Encyclopedia of Mathematics and its Applications}.
\newblock Cambridge University Press, Cambridge, 1993.

\bibitem[Sch07]{jms}
J.-M. Schlenker.
\newblock Small deformations of polygons and polyhedra.
\newblock {\em Trans. Amer. Math. Soc.}, 359(5):2155--2189 (electronic), 2007.

\bibitem[Thu97]{Tlivre}
W.~P. Thurston.
\newblock {\em Three-dimensional geometry and topology. {V}ol. 1}, volume~35 of
  {\em Princeton Mathematical Series}.
\newblock Princeton University Press, Princeton, NJ, 1997.
\newblock Edited by Silvio Levy.

\bibitem[Thu98]{T}
W.~P. Thurston.
\newblock Shapes of polyhedra and triangulations of the sphere.
\newblock In {\em The Epstein birthday schrift}, volume~1 of {\em Geom. Topol.
  Monogr.}, pages 511--549 (electronic). Geom. Topol. Publ., Coventry, 1998.
\newblock Circulated as a preprint since 1987.

\bibitem[Tro86]{troyanov2}
M.~Troyanov.
\newblock Les surfaces euclidiennes \`a singularit\'es coniques.
\newblock {\em Enseign. Math. (2)}, 32(1-2):79--94, 1986.

\bibitem[Tro91]{troyanov1}
M.~Troyanov.
\newblock Prescribing curvature on compact surfaces with conical singularities.
\newblock {\em Trans. Amer. Math. Soc.}, 324(2):793--821, 1991.

\bibitem[Tro07]{Troteich}
M.~Troyanov.
\newblock On the moduli space of singular {E}uclidean surfaces.
\newblock In {\em Handbook of Teichm\"uller theory. Vol. I}, volume~11 of {\em
  IRMA Lect. Math. Theor. Phys.}, pages 507--540. Eur. Math. Soc., Z\"urich,
  2007.

\bibitem[Tum07]{tumarkin}
P.~Tumarkin.
\newblock Compact hyperbolic {C}oxeter {$n$}-polytopes with {$n+3$} facets.
\newblock {\em Electron. J. Combin.}, 14(1):Research Paper 69, 36 pp.
  (electronic), 2007.

\bibitem[Vin75]{vinberg2}
{\`E}.~B. Vinberg.
\newblock Some arithmetical discrete groups in {L}oba\v cevski\u\i\ spaces.
\newblock In {\em Discrete subgroups of Lie groups and applications to moduli
  (Internat. Colloq., Bombay, 1973)}, pages 323--348. Oxford Univ. Press,
  Bombay, 1975.

\bibitem[Vin85]{vinberg}
{\`E}.~B. Vinberg.
\newblock Hyperbolic groups of reflections.
\newblock {\em Russian Math. Surveys}, 40(1):31--75, 1985.

\bibitem[vin93]{vinberglivre}
{\em Geometry. {II}}, volume~29 of {\em Encyclopaedia of Mathematical
  Sciences}.
\newblock Springer-Verlag, Berlin, 1993.
\newblock Spaces of constant curvature, A translation of Geometriya. II, Akad.
  Nauk SSSR, Vsesoyuz. Inst. Nauchn. i Tekhn. Inform., Moscow, 1988,
  Translation by V. Minachin [V. V. Minakhin], Translation edited by \`E. B.
  Vinberg.

\bibitem[Web93]{weber}
M.~Weber.
\newblock {\em Fundamentalbereiche komplex hyperbolischer {F}l\"achen}.
\newblock Bonner Mathematische Schriften [Bonn Mathematical Publications], 254.
  Universit\"at Bonn Mathematisches Institut, Bonn, 1993.

\bibitem[YNK02]{kojimaal2}
Y.~Yamashita, H.~Nishi, and S.~Kojima.
\newblock Configuration spaces of points on the circle and hyperbolic {D}ehn
  fillings. {II}.
\newblock {\em Geom. Dedicata}, 89:143--157, 2002.

\bibitem[Yos96]{yoshida1}
M.~Yoshida.
\newblock The democratic compactification of configuration spaces of point sets
  on the real projective line.
\newblock {\em Kyushu J. Math.}, 50(2):493--512, 1996.

\bibitem[Yos01]{yoshida2}
M.~Yoshida.
\newblock A hyperbolic structure on the real locus of the moduli space of
  marked cubic surfaces.
\newblock {\em Topology}, 40(3):469--473, 2001.

\end{thebibliography}

\end{document}